\documentclass [a4paper, 10 pt]{article}
\usepackage{amsmath}
\usepackage{amsfonts}
\usepackage{amsthm}
\usepackage{amscd}
\usepackage{amssymb}
\usepackage{graphicx}
\usepackage{epsfig}
\usepackage{graphics}

\newcommand\I{{\bf I}}
\newcommand {\Ical}{{\cal{{I}_{LG}}}}
\newcommand\grad{{\bf \nabla}}
\newcommand\rvec{{\bf r}}
\newcommand\la{{\lambda}}
\newcommand\ga{{\gamma}}
\newcommand\n{{\bf n}}
\newcommand\nvec{{\bf n}}
\newcommand\x{{\bf x}}
\newcommand\uvec{{\bf u}}
\newcommand\zvec{{\bf z}}
\newcommand\xvec{{\bf x}}

\newcommand\mvec{{\bf m}}
\newcommand\evec{{\bf e}}

\newcommand\Pvec{{\bf P}}

\newcommand\eps{{\epsilon}}
\newcommand\Q{{\bf Q}}
\newcommand\Qvec{{\bf Q}}
\newcommand{\Acal}{{\cal A}}

\newcommand{\Rr}{{\mathbb R}}

\newtheorem{lem}{Lemma}[section]
\newtheorem{prop}{Proposition}[section]

\theoremstyle{definition}

\theoremstyle{remark}

\def\comment#1{}

\def\withcomments{
\addtolength{\oddsidemargin}{-0.5 in}
\addtolength{\evensidemargin}{-0.5 in}
\newcounter{mycommentcounter}
\def\comment##1{\refstepcounter{mycommentcounter}%
  \ifhmode%
  \unskip%
  {\dimen1=\baselineskip \divide\dimen1 by 2 %
    \raise\dimen1\llap{\tiny -\themycommentcounter-}}\fi%
  \marginpar{\renewcommand{\baselinestretch}{0.8}%
    \footnotesize [\themycommentcounter]: \raggedright ##1}}
}

\withcomments

\begin{document}

\title{The Radial-Hedgehog Solution in Landau--de Gennes' theory}
 \date{\today}
\author{Apala Majumdar}
 \maketitle

\begin{abstract} We study the radial-hedgehog solution
in a three-dimensional spherical droplet, with homeotropic
boundary conditions, within the Landau-de Gennes theory for
nematic liquid crystals. The radial-hedgehog solution is a
candidate for a globally stable configuration in this model
framework and is also a prototype configuration for studying
isolated point defects in condensed matter physics. The static
properties of the radial-hedgehog solution are governed by a
nonlinear singular ordinary differential equation. We consider two
different limits separately - the \emph{vanishing core limit} and
\emph{low-temperature limit} respectively. We use a combination of
Ginzburg-Landau techniques, perturbation methods and stability
analysis to study the qualitative properties of the
radial-hedgehog solution, both in the vicinity of and away from
the defect core. We establish the instability of the
radial-hedgehog solution with respect to biaxial perturbations in
certain parameter regimes and demonstrate the stability of the
radial-hedgehog solution in other parameter regimes. Our results
complement previous work in the field, are rigorous in nature and
give information about the role of geometry and temperature on the
properties of the radial-hedgehog solution and the associated
biaxial instabilities.
\end{abstract}

\section{Introduction}
\label{sec:intro}

Defect structures have attracted a lot of interest in the liquid
crystal community \cite{penzenstadler,rossovirga,schopohl,sonnet}.
Defect structures in liquid crystalline systems are usually
modelled within the Landau-de Gennes framework, whereby the liquid
crystal configuration is mathematically described by a symmetric,
traceless $3\times 3$ matrix, known as the $\Q$-tensor order
parameter \cite{dg}. The $\Q$-tensor can be written in terms of
its eigenvalues and eigenvectors as shown below -
\begin{equation}
\label{eq:int1}
\Q = \sum_{i=1}^{3}\la_i \evec_i\otimes \evec_i, \quad \sum_i \la_i = 0
\end{equation}
where $\la_i$ are the eigenvalues and $\evec_i$ are the
corresponding orthonormal eigenvectors. The liquid crystal is said
to be in the (i) isotropic state when $\la_i=0$ for $i=1\ldots 3$,
(ii) uniaxial state when $\Q$ has a pair of equal non-zero
eigenvalues and (iii) biaxial state when $\Q$ has three distinct
non-zero eigenvalues \cite{newtonmottram}.

A prototype example of such a confined system is a spherical
droplet with \emph{strong radial anchoring} or \emph{homeotropic}
(normal) boundary conditions. This example has been widely studied
in the literature, especially from a numerical point of view, and
it is generally believed that there are two competing equilibrium
configurations - (a)the \emph{radial-hedgehog} solution which has
a single isolated point defect at the droplet centre and (b) the
\emph{biaxial-torus} solution where the point defect broadens out
to a ring-like structure around the droplet centre
\cite{mkaddem&gartland2,kralj,mkaddem&gartland1,schopohl,sonnet}.
The radial-hedgehog solution is purely uniaxial everywhere except
for an isotropic point at the droplet centre whereas the
biaxial-torus configuration exhibits a high degree of biaxiality
around the droplet centre. The isotropic point in the
radial-hedgehog solution and the biaxial ring in the torus
solution are interpreted as being \emph{defect structures} since
they are localised regions of abrupt changes in the eigenvalue
structure.

This paper aims to build a self-contained mathematical description
of the radial-hedgehog solution within the Landau-de Gennes
framework. Firstly, this is an interesting mathematical problem in
its own right since the radial-hedgehog solution is a rare example
of an explicit solution of the Landau-de Gennes Euler-Lagrange
equations in (\ref{eq:7}). Moreover, the corresponding scalar
order parameter is a solution of an ordinary differential equation
(see (\ref{eq:11})) and hence, has a tractable and yet non-trivial
mathematical structure. Indeed, this is the first step in the
mathematical theory of defects in liquid crystalline systems.
Secondly, a systematic mathematical analysis of the
radial-hedgehog solution is crucial for understanding the
structure and locations of point defects in liquid crystalline
systems, the multiplicity of uniaxial solutions and the
characterization of the competing biaxial structures.

To further elaborate on the above, radial-hedgehog solutions can
be thought of as prototypical vortices in the Ginzburg-Landau
theory for superconductors \cite{bbh2}. More precisely, the
radial-hedgehog solution can be interpreted as being a degree $+1$
vortex in three dimensions. There is a very well-developed theory
for the structure, location, multiplicity and stability of
vortices in Ginzburg-Landau theory, especially in two dimensions
but generalizations to higher dimensions are non-trivial
\cite{bbh2,millot,farina,gustafson}. One of the main aims of this
paper is to clearly demonstrate the analogies and differences
between the mathematical formulation of radial-hedgehog solutions
in the Landau-de Gennes framework and Ginzburg-Landau vortices.
Once the inter-relationship is correctly understood, this will
contribute to a sound theoretical foundation for defects in liquid
crystals and Ginzburg-Landau numerical methods can also be used
for the simulation of defects in liquid crystalline systems. We
deal with two separate limits in this paper - the
\emph{low-temperature limit} where the governing ordinary
differential equation has an \emph{almost} Ginzburg-Landau
structure and the \emph{vanishing core limit} where there are
important technical differences. In particular, we cannot exploit
Ginzburg-Landau techniques to describe the isotropic defect core
in the vanishing core limit. The low-temperature limit is relevant
for liquid crystalline systems \emph{deep} in the nematic phase
where we expect to see a high degree of orientational ordering.
The vanishing core limit is relevant for materials whose elastic
constants are typically much smaller in magnitude than the
thermotropic parameters and quoted values in the literature
suggest that this limit is relevant for commonly used liquid
crystalline materials \cite{elastic}. More generally, although the
study of uniaxial states can be viewed as a generalized
Ginzburg-Landau theory from $\Rr^3 \to \Rr^3$ \cite{maj1},
biaxiality presents a whole host of new mathematical challenges,
outside the scope of Ginzburg-Landau theory \cite{amaz}. In
particular, there is no analogue of a biaxial instability in the
current Ginzburg-Landau literature and such instabilities play a
pivotal role in Landau-de Gennes theory.

The paper is organized as follows. In Section~\ref{sec:prelim}, we
prove the existence of a radial-hedgehog solution in the Landau-de
Gennes framework and establish bounds for the corresponding scalar
order parameter. In Section~\ref{sec:isotropic}, we derive a
series expansion for the radial-hedgehog solution near its
isotropic core and demonstrate its similarity with
three-dimensional vortices in Ginzburg-Landau theory
\cite{farina}. We then show that the radial-hedgehog solution
cannot be a global Landau-de Gennes energy minimizer for
sufficiently large droplets, for sufficiently low temperatures by
means of an explicit comparison argument . This result is
qualitatively similar to a result reported in
\cite{mkaddem&gartland2} but our method of proof is different. In
Section~\ref{sec:lowtemp}, we focus on the low-temperature limit
and the resulting Ginzburg-Landau structure of the governing
ordinary differential equation. We use shooting arguments to
establish qualitative properties of the corresponding scalar order
parameter and use Ginzburg-Landau techniques to prove the
uniqueness of the radial-hedgehog solution in this limit. We deal
with the vanishing core limit in Section~\ref{sec:ns}, whereby the
governing ordinary differential equation does not have a
Ginzburg-Landau structure. Our results are weaker in this case and
describe the far-field properties, away from the defect core. In
Section~\ref{sec:stability}, we perform a linear stability
analysis of the radial-hedgehog solution and this stability
analysis gives insight into the effect of the ball radius on the
associated equilibrium structure. In Section~\ref{sec:discussion},
we discuss our results and how they complement previous work in
this area.

\section{Preliminaries}
\label{sec:prelim}

We study the qualitative properties of radial-hedgehog solutions
inside spherical droplets, $B(0,R)\subset \Rr^3$, where
\begin{equation}
\label{eq:droplet}
B(0,R) = \left\{\rvec\in\mathbb{R}^3; |\rvec|\leq R \right\}
\end{equation}
and $R>0$ is independent of any model parameters, subject to
strong radial anchoring conditions. We work within the Landau-de
Gennes theory for nematic liquid crystals, in the low-temperature
regime.

Let $\bar{S}\subset \mathbb{M}^{3\times 3}$ denote the space of
symmetric, traceless $3\times 3$ matrices  i.e.
\begin{displaymath}
\label{eq:4a}
 \bar{S}\stackrel{def}{=} \left\{\Q \in \mathbb{M}^{3\times 3};
\Q_{ij}=\Q_{ji},~\Q_{ii} = 0 \right\}
\end{displaymath}
where we have used the Einstein summation convention; the Einstein
convention will be used in the rest of the paper. The
corresponding matrix norm is defined to be
\begin{displaymath}
\label{eq:4b}
 \left| \Q \right|\stackrel{def}{=}\sqrt{\textrm{tr}\Q^2} =\sqrt{ \Q_{ij}
\Q_{ij}} \quad i,j=1\ldots 3.
\end{displaymath}We recall from \cite{amaz,newtonmottram} that
an arbitrary $\Q\in\bar{S}$ can be written as
$$\Q =s\left(\n\otimes\n - \frac{1}{3}\mathbf{I}\right) + r \left(\mvec\otimes \mvec - \frac{1}{3}\mathbf{I}\right)$$
where $\n,\mvec$ are orthonormal eigenvectors of $\Q$,
$s,r$ are scalar order parameters and we either have $0\leq r\leq \frac{s}{2}$ or $\frac{s}{2}\leq r \leq 0$. If $\Q\in\bar{S}$ is uniaxial, then this representation formula can be simplified to
$$\Q = s\left( \n\otimes\n - \frac{1}{3}\mathbf{I}\right)$$
where $\n$ is the leading eigenvector of $\Q$ and $s$ is
a scalar order parameter that measures the degree of orientational ordering about $\n$.

The Landau-de Gennes energy functional is given by \cite{dg,newtonmottram}
\begin{equation}
\label{eq:1}
\Ical[\Q] = \int_{B(0,R)} \frac{L}{2}\left|\grad \Q\right|^2 + f_B(\Q)~dV
\end{equation} where $\left|\grad \Q\right|^2 = \sum_{i,j,k=1}^{3}\left(\frac{\partial\Q_{ij}}{\partial \x_k}\right)^2$ is the elastic energy density, $L$ is a small material-dependent elastic constant and $f_B$ is the bulk energy density given by
\begin{equation}
\label{eq:2}
f_B(\Q) = -\frac{a^2}{2}\textrm{tr}\Q^2 - \frac{b^2}{3}\textrm{tr}\Q^3 + \frac{c^2}{4}\left(\textrm{tr} \Q^2 \right)^2 + C\left(a^2,b^2,c^2\right).
\end{equation} The form (\ref{eq:2}) is the simplest
form of the bulk energy density that allows for a first-order
nematic-isotropic phase transition; here $b^2,c^2$ are
material-dependent positive constants and $a^2>0$ is a
temperature-dependent parameter. Typical values of these
characteristic constants are $a^2 = 0.042\times 10^6 (T^* - T)N/K
m^2, b^2 = 0.64\times 10^6 N/m^2, c^2 = 0.35\times 10^6 N/m^2$
where $T$ is the absolute temperature and $T^*$ is a
characteristic temperature below which the isotropic phase
$\Qvec=0$ ceases to be a locally stable stationary point of
$f_B$ in (\ref{eq:2}) \cite{newtonmottram, elastic}. 
We work in the low-temperature regime where
the bulk energy density attains its global minimum on the set of
uniaxial $\Q$-tensors given by \cite{maj1}
\begin{eqnarray}
\label{eq:3} && \Q_{min} =
 \left\{\Q\in \bar{S}, ~ \Q= s_+ \left( \n\otimes \n -
\frac{1}{3}\I \right)~
\right\}\end{eqnarray} with $\n\in\mathbb{S}^2$ and
\begin{equation}
s_+=\frac{b^2+\sqrt{b^4+24a^2c^2}}{4c^2}.
\label{eq:4}
\end{equation} In particular, as $a^2$ increases, we move to
lower temperatures \textit{deep} in the nematic phase.
The additive constant $C\left(a^2,b^2,c^2\right)$ in (\ref{eq:2}) ensures that
$f_B(\Q)\geq 0$ for all $\Q\in\bar{S}$. The admissible space is defined to be
\begin{eqnarray}
&& \Acal_{\Q} = \left\{\Q\in W^{1,2}\left(B(0,R);\bar{S}\right);
\textrm{$\Q=\Q_b$ on $\partial B(0,R)$} \right\}\label{eq:5}
\end{eqnarray}  and the Dirichlet boundary condition $\Q_b \in \Q_{\min}$ is specified to be
\begin{equation}
\label{eq:6} \Q_b = s_+\left(\frac{\rvec}{|\rvec|} \otimes
\frac{\rvec}{|\rvec|} - \frac{1}{3}\I\right).
\end{equation} This is referred to as \emph{strong radial anchoring}
in the liquid crystal literature
\cite{mkaddem&gartland1,mkaddem&gartland2}, since
$\frac{\rvec}{|\rvec|}$ is the unit vector in the radial
direction. The physically observable, equilibrium configurations
correspond to either global or local minimizers of $\Ical$ in
$\Acal_\Q$. For completeness, we recall that
$W^{1,2}\left(B(0,R);\bar{S}\right)$ is the Sobolev space of
square-integrable $\Q$-tensors with square-integrable first
derivatives \cite{evans}. The corresponding $W^{1,2}$-norm is
given by $\| \Q \|_{W^{1,2}(B(0,R))} =\left( \int_{B(0,R)} |\Q|^2
+ |\grad \Q|^2~dx\right)^{1/2}.$ In addition to the
$W^{1,2}$-norm, we also use the $L^{\infty}$-norm in this paper,
defined to be $\|\Q\|_{L^{\infty}(B(0,R))} = \textrm{ess
sup}_{\x\in B(0,R)}|\Q(\x)|$.

In what follows, we consider two different limits: the $L\to 0$
limit which is referred to as the \emph{vanishing core limit} and
the $a^2 \to \infty$ limit, which is referred to as the
\emph{low-temperature limit}. The reason for making this
distinction will become clear in the subsequent sections. We work
in a dimensionless framework and as outlined in
\cite{mkaddem&gartland1,mkaddem&gartland2}, we introduce the
following dimensionless variables : -
\begin{eqnarray}
\label{eq:nondim}
&& \bar{\rvec} = \frac{\rvec}{\xi}, \quad \bar{\Qvec} = \sqrt{\frac{27 c^4}{ 2b^4}} \Qvec, ~ \overline{\Ical} = \left(\frac{27 c^6}{4 b^4 L^3}\right)\Ical
\end{eqnarray}
where $\xi = \sqrt{\frac{27c^2 L}{b^4}}$. It is straightforward to
show that the corresponding dimensionless energy density is
\begin{eqnarray}
\label{eq:nondim2} \bar{e}(\bar{\Qvec},\grad\bar{\Qvec}) =
\frac{1}{2}|\grad \bar{\Qvec}|^2 -
\frac{t}{2}\textrm{tr}\bar{\Qvec}^2 -
\sqrt{6}\textrm{tr}\bar{\Qvec}^3 +
\frac{1}{2}\left(\textrm{tr}\bar{\Qvec}^2\right)^2 + C(t)
\end{eqnarray}
where $t = \frac{27 a^2 c^2}{b^4}>0$ is a dimensionless
\emph{reduced temperature}, $C(t)$ is an additive constant that
ensures $f_B(t)=- \frac{t}{2}\textrm{tr}\bar{\Qvec}^2 -
\sqrt{6}\textrm{tr}\bar{\Qvec}^3 +
\frac{1}{2}\left(\textrm{tr}\bar{\Qvec}^2\right)^2 + C(t)\geq 0$
and $t>1$ throughout the paper.

We employ a second change of variable
\begin{eqnarray}\label{eq:nondimnew}
\tilde{\Qvec} = \frac{\bar{\Qvec}}{h_+}; \quad \tilde{\rvec} = \sqrt{t}\bar{\rvec}
\end{eqnarray}
where
\begin{eqnarray}
\label{eq:h+}
h_+ = \frac{ 3 + \sqrt{9 + 8t}}{4}.
\end{eqnarray}
This choice of the dimensionless variables is especially relevant
for the low-temperature limit, as will be demonstrated in
Section~\ref{sec:lowtemp}. The corresponding dimensionless energy
density is
\begin{eqnarray}
\label{eq:nondim3}
\tilde{e}(\tilde{\Qvec},\grad\tilde{\Qvec}) = \frac{1}{2}|\grad\tilde{\Qvec}|^2 - \frac{1}{2}\textrm{tr}\tilde{\Qvec}^2 - \frac{\sqrt{6}h_+}{t}\textrm{tr}\tilde{\Qvec}^3 + \frac{h_+^2}{2t}\left(\textrm{tr}\tilde{\Qvec}^2\right)^2.
\end{eqnarray}
One can readily compute the Euler-Lagrange equations associated with the energy functional,
\begin{equation}
\label{eq:nondim4}
\tilde{\mathcal{I}}_{LG}[\tilde{\Qvec}] = \int_{B(0,\tilde{R})}\tilde{e}(\tilde{\Qvec},\grad\tilde{\Qvec})~dV,
\end{equation}
where $\tilde{R} = \sqrt{t}\frac{R}{\xi}$. In what follows, we drop the
 \emph{tilde} on the dimensionless variables for brevity and all subsequent results are to be understood
 in terms of the dimensionless variables.
The associated Euler-Lagrange equations are \cite{maj1,amaz} -
\begin{equation}
\Delta \Q_{ij}=  - \Q_{ij}  - \frac{ 3\sqrt{6} h_+}{t}\left(\Q_{ik}\Q_{kj}-\frac{\delta_{ij}}{3}\textrm{tr}(\Q^2)\right)
+\frac{2 h_+^2}{t}\Q_{ij}\textrm{tr}(\Q^2),~ ~i,j=1,2,3 \label{eq:7},
\end{equation} where the term $\frac{\delta_{ij}}{3}\textrm{tr}(\Q^2)$ is a Lagrange multiplier
associated with the tracelessness constraint. It follows from
standard arguments in elliptic regularity that any solution $\Q$
of the nonlinear elliptic system (\ref{eq:7}) is smooth and real
analytic on $B(0,R)$ \cite{amaz}. In
particular, all global and local energy minimizers in $\Acal_\Q$
are classical solutions of (\ref{eq:7}).

Radial-hedgehog solutions are examples of spherically-symmetric
uniaxial solutions of the system (\ref{eq:7}) in the admissible
space $\Acal_{\Q}$ and have the form
\begin{equation}
\label{eq:8} \Q = \sqrt{\frac{3}{2}} h(r)
\left(\frac{\rvec}{|\rvec|} \otimes \frac{\rvec}{|\rvec|} -
\frac{1}{3}\I\right).
\end{equation}
Here the scalar order parameter $h$ only depends on the radial
distance $r = |\rvec|$ from the origin and the corresponding
admissible space is defined to be \begin{equation} \label{eq:9}
\Acal_h = \left\{ h\in W^{1,2}\left([0,R], \Rr\right); h(R) =
1\right\}.\end{equation} We note that $\Q \in
W^{1,2}\left(B(0,R);\bar{S}\right)$ necessarily implies that $h\in
W^{1,2}\left([0,R];\Rr\right)$ since the eigenvalues of a
symmetric matrix are Lipschitz functions of the matrix components
\cite{strongsemismoothness} and hence, $\Acal_h$ is a natural
choice for the admissible space. There may be multiple
spherically-symmetric solutions of (\ref{eq:7}) but we define a
radial-hedgehog solution to be an energy-minimizing
spherically-symmetric solution in the sense described below.

\begin{prop}
\label{prop:1}
(a) Consider the energy functional
\begin{equation}
\label{eq:10} I[h] = \int_{0}^{R} r^2 \left( \frac{1}{2}
\left(\frac{dh}{dr}\right)^2 + \frac{3 h^2}{r^2} + f(h) \right)~dr
\end{equation} defined for functions $h \in \Acal_h$, where
\begin{equation}
\label{eq:fh} f\left(h \right) = - \frac{h^2}{2} -
\frac{h_+}{t}h^3 + \frac{h_+^2}{2t}h^4 + C(t), \end{equation}
$f(h)\geq 0$ for $h\in\Acal_h$ and $f(1)=0$. There \emph{exists} a
global minimizer $h^* \in \Acal_h$ for $I$ in (\ref{eq:10}). The
function $h^*$ is a solution of the following singular nonlinear
ordinary differential equation
\begin{equation}
\label{eq:11}
\frac{d^2 h}{dr^2} + \frac{2}{r}\frac{dh}{dr} - \frac{6h}{r^2} = -h + h^3 + \frac{3 h_+}{t}\left(h^3 - h^2\right)
\end{equation}
subject to the boundary conditions
\begin{equation}
\label{eq:12}
h(0) = 0 ~and~ h(R) = 1.
\end{equation} Moreover, $h^*$ is analytic for all $r\geq 0$.

\vspace{.1 cm}
(b) Define the radial-hedgehog solution by
\begin{equation}
\label{eq:13}
\Q^* = \sqrt{\frac{3}{2}} h^*(r)\left(\frac{\rvec\otimes \rvec}{r^2} - \frac{1}{3}\I\right)
\end{equation} where $h^*$ is a global minimizer of $I[h]$ in (\ref{eq:10}),
in the admissible space $\Acal_h$ . Then $\Q^*$ is a solution of
the Landau-de Gennes Euler-Lagrange equations (\ref{eq:7}) i.e. is
a stationary point of the Landau-de Gennes energy functional.
Moreover, these solutions satisfy the following energy bound
\begin{equation}
\label{eq:14}
\tilde{\mathcal{I}}_{LG}[\Q^*]
\leq 12 \pi  R 
\end{equation} where $\tilde{\mathcal{I}}$ has been defined in (\ref{eq:nondim4}).

(d) The function $h^*$ satisfies the following bounds for $r\in[0,R]$ -
\begin{equation}
\label{eq:15}
0\leq h^*(r) \leq 1 \quad r\in[0,R].
\end{equation}
\end{prop}

\textit{Proof}

(a) Consider the energy functional $I[h]$ defined for
$h\in\Acal_h$. Firstly, we note that the admissible space
$\Acal_h$ is non-empty. Indeed, the constant function $h(r)=1$ for
$r\in [0,R]$ belongs to $\Acal_h$. Secondly, the functional $I$ in
(\ref{eq:10}) is bounded from below i.e. $I[h]\geq 0$ for $h\in
\Acal_h$ and is weakly lower semicontinuous on our admissible
space (since the integrand is convex in $dh/dr$). The existence of
a global minimizer $h^* \in\Acal_h$ now follows from the direct
methods in the calculus of variations \cite{evans}.

It is straightforward to compute the Euler-Lagrange equations associated with the functional $I$ in (\ref{eq:10}) i.e.
$$ \frac{d}{dr} \left(\frac{\partial e(h, h^{'})}{\partial h'}\right) = \frac{\partial e(h, h')}{\partial h} $$
where $h' = dh/dr$, $e(h, h') = r^2 \left( \frac{1}{2} \left(\frac{dh}{dr}\right)^2 + \frac{3 h^2}{r^2} - \frac{h^2}{2} - \frac{h_+}{t}h^3 + \frac{h_+^2}{2t}h^4 \right)$. One can check that the corresponding Euler-Lagrange equation is indeed the ordinary differential equation in (\ref{eq:11}) and a global minimizer $h^*$ is necessarily a solution of these Euler-Lagrange equations.

The boundary condition $h^*(R) = 1$ follows from our definition of
the admissible space $\Acal_h$. All functions $h\in \Acal_h$ are
necessarily continuous since $h\in W^{1,2}([0,R],\Rr) \implies
h\in C^{0,\alpha}([0,R],\Rr)$ for some $0<\alpha<1/2$ from the
Sobolev embedding theorem \cite{evans}. The boundary condition
$h(0)=0$ follows from the continuity of $h^*(r)$ for $r\in[0,R]$.
We proceed by contradiction and assume that  $|h^*(r)|\geq h_0$
for $r\in [0,r_0]$, for some fixed $h_0>0$ and $0<r_0<<1$. Since
$h^*$ is continuous, we deduce that $h^*$ has a fixed sign near
the origin and we further assume that $h^*(r)>h_0>0$ for
$r\in[0,r_0]$. Consider the governing equation (\ref{eq:11}); it
can be re-written as
\begin{equation}
\frac{d}{dr}\left( r^2 \frac{dh}{dr}\right) = 6h + r^2\left( h^3 - h + \frac{3 h_+}{t}\left(h^3 - h^2\right)\right)
\label{eq:16}
\end{equation} where $h_+$ has been defined in (\ref{eq:h+}).
Then, we have
\begin{equation}
\label{eq:17}
r^2 \frac{dh}{dr} \geq \int_{\eps}^{r} 6 h(r') dr' + C r^3 + \eps^2 h^{'}(\eps) \quad ~for~ r\in(0,r_0)
\end{equation}
where $0<\eps<r/10$ is fixed, $h^{'}(\eps) = \frac{dh}{dr}|_{r=\eps}$ and $C$ is a constant independent of $\eps$. We note that $h^{'}(\eps)$ can be bounded independently of $\eps$ i.e. $\left|\frac{dh}{dr}\right|\leq C(t)$ for $r\in\left[0,R\right]$ from \cite{amaz}. Squaring both sides of (\ref{eq:17}) and integrating from $\eps$ to $r$, we obtain
\begin{equation}
\label{eq:18}
 \int_{\eps}^{r} \left(\frac{dh}{dr'}\right)^2 ~dr'
 \geq \int_{\eps}^{r} \frac{\ga h_0^2}{t^2} ~dt + C^{''} r^3 +\eps^2 h^{'}(\eps) \int_{\eps}^{r}\frac{1}{t^3}dt \quad for~r\in(0,r_0),
\end{equation} where $C^{''}$ is a constant independent of $\eps$.
In the limit $\eps \to 0$, (\ref{eq:18}) contradicts the hypothesis that $h\in W^{1,2}\left([0,R];\Rr\right)$ from which we must have
$$  \int_{0}^{R}\left(\frac{dh}{dr}\right)^2 ~dr < \infty.$$
Therefore, we deduce that $h(0)=0$ for any solution of
(\ref{eq:11}) in $\Acal_h$ and $h^*\in\Acal_h$ is a solution of
(\ref{eq:11}), subject to the boundary conditions (\ref{eq:12}).
The analyticity of $h^*$ now follows from standard arguments in
the theory of ordinary differential equations \cite{herve&herve}.

(b) Given $h^*$, define a radial-hedgehog solution as follows
$$\Q^* = \sqrt{\frac{3}{2}}
h^*(r)\left(\frac{\rvec\otimes\rvec}{r^2} -
\frac{1}{3}\mathbf{I}\right).$$ From \cite{strongsemismoothness}, $\Q^*\in W^{1,2}(B(0,R);\bar{S})$ necessarily implies $h^*\in W^{1,2}([0,R];\Rr^+)$ and the preceding arguments necessarily require $h^*(0)=0$. One can directly check that
\begin{equation}
\label{eq:energybound2}
\Ical[\Q^*] = 4 \pi I[h^*]
\end{equation}
and that $\Q^*$ is a solution of the Euler-Lagrange equations (\ref{eq:7}), since $h^*$ is a solution of the ordinary differential equation (\ref{eq:11}), subject to the boundary conditions (\ref{eq:12}).

The function $h^*$ has been defined to be the global minimizer of
the functional $I$ in (\ref{eq:10}), in the admissible space
$\Acal_h$. However, the constant function, $\bar{h}(r) = 1$ for
$r\in[0,R]$, belongs to $\Acal_h$ and hence
\begin{equation}\label{eq:19}
I[h^*] \leq I[\bar{h}] =3R .
\end{equation} The energy bound on $\Ical[\Q^*]$, where $\Q^* = \sqrt{\frac{3}{2}}h^*(r)\left(\frac{\rvec\otimes\rvec}{r^2} - \frac{1}{3}\mathbf{I}\right)$, follows from (\ref{eq:energybound2}).

(c) The upper bound $|h^*(r)|\leq 1$ follows directly from a
result in \cite{maj2} where we establish that every solution $\Q$
of the system (\ref{eq:7}) in the admissible space $\Acal_\Q$
satisfies the global upper bound
$$ \left| \Q (\rvec)\right| \leq1.$$ The radial-hedgehog solution $\Q^*$ is a solution of the system (\ref{eq:7}) and
$$ \left| \Q^*(\rvec)\right| =  |h^*(r)|$$
where $r=\left|\rvec\right|$. The upper bound $|h^*(r)|\leq 1$ follows immediately.

The lower bound $h^*(r)\geq 0$ follows from the energy minimality condition. We assume that there exists an interior measurable subset
$$\tilde{\Omega} = \left\{\rvec\in\Omega;~ h^*(r) < 0\right\} \subset B(0,R)$$ with $h^*(r) = 0$ on $\partial\tilde{\Omega}$. We note that $\tilde{\Omega}$ must be an interior subset because of the boundary condition $\Q_b$ in (\ref{eq:6}). We define the perturbation
\begin{equation} \label{eq:20}
 \bar{h}^*(r) =
  \begin{cases}
   h^*(r), & \rvec\in B(0,R)\setminus \tilde{\Omega},\\
  - h^*(r),&
    \rvec\in\tilde{\Omega}.
  \end{cases}
\end{equation} One can then easily check that
\begin{equation}\label{eq:21}
I[\bar{h}^*] - I[h^*] = \int_{\tilde{\Omega}} -\frac{h_+}{3}
\left(\bar{h}^*\right)^3 +\frac{h_+}{3} h{^*}^3~dV =
\int_{\tilde{\Omega}} \frac{2 h_+}{t}h{^{*}}^3~dV <0
\end{equation}
since $h^*(r)<0$ on $\tilde{\Omega}$ by assumption. The inequality
(\ref{eq:21}) contradicts the global minimality of $h^*$ in
$\Acal_h$ and hence, we deduce that $h^*(r)\geq 0$ for
$r\in[0,R]$. The inequalities (\ref{eq:15}) now follow. $\Box$

In summary, in Proposition~\ref{prop:1}, we prove the existence of
a radial-hedgehog solution of the form (\ref{eq:13}), that can be
interpreted as being a Landau-de Gennes energy minimizer within
the class of radially-symmetric configurations. This
radial-hedgehog solution satisfies the energy bound (\ref{eq:14})
and the corresponding scalar order parameter $h^*$ is bounded from
both above and below as shown in (\ref{eq:15}).  The
radial-hedgehog solution has a single isolated isotropic point at
the origin where $h^*$ vanishes and this isolated isotropic point
is interepreted as being a defect point, since the radial-hedgehog
solution is strictly uniaxial everywhere else. In the next
section, we study the \emph{isotropic core} of the radial-hedgehog
solution and the manifestation of biaxial instabilities within
this core.

\section{The Isotropic Core and Biaxial Instabilities}
\label{sec:isotropic}

\begin{prop}\label{prop:iso2} Let $h^*$ be a global minimizer of
the energy functional $I$ in (\ref{eq:10}). Then $h^*$ is a
solution of the ordinary differential equation (\ref{eq:11})
subject to the boundary conditions (\ref{eq:12}). As $r \to 0$, we
have the following series expansion for $h^*$
\begin{equation}
\label{eq:iso2}
h^*(r)= \sum_{n=0}^{\infty} a_n r^n = a_2 r^2 \left[1 - \frac{r^2}{14}+ o(r^2)\right] \quad as~r\to 0
\end{equation}
where $a_n = 0$ for all $n$ odd and $a_2>0$ is an arbitrary constant.
In addition, as the domain size $R \to \infty$, we also have that \cite{mkaddem&gartland1,mkaddem&gartland2}
\begin{equation}
\label{eq:bounds}
\frac{r^2}{r^2 + 14} \leq h^*(r) \leq \frac{r^2}{r^2 + t \lambda_t^2}
\end{equation} where $\lambda_t^2 = \frac{24}{9 + 8t + 3\sqrt{9+8t}} \leq \frac{3}{t}\leq 3$
since $t\geq 1$.
Therefore, for $R$ sufficiently large, we have the following bounds on the constant $a_2$ in (\ref{eq:iso2})
\begin{equation}
\label{eq:bounds2}
 \frac{1}{14} \leq a_2 \leq \frac{1}{3} + \frac{3}{8t} + \frac{1}{8t}\sqrt{9+8t}.
 \end{equation}
\end{prop}
\textit{Comment: Equation~(\ref{eq:iso2}) is identical to the
series expansion for three-dimensional vortices near the origin,
within the Ginzburg-Landau theory for superconductivity
\cite{farina}.}

\textit{Comment: The limit $R\to \infty$ corresponds to either the
$L\to 0$ limit in (\ref{eq:1}) or the $t\to \infty$ limit in
(\ref{eq:nondim3}) or the doubly infinite limit $L\to 0, t\to
\infty$.}

\textit{Proof:} From Proposition~\ref{prop:1}, we have that $h^*$
is analytic for $r\geq0$. We seek a power series expansion of
$h^*$ around the origin with $h^*(0)=0$, of the form
\begin{equation}
\label{eq:iso6}
h^*(r) = \sum_{n=1}^{\infty} a_n r^n \quad 0< r \leq R_0
\end{equation} where $R_0$ is the radius of convergence.

We substitute the ansatz (\ref{eq:iso6}) into the ordinary
differential equation (\ref{eq:11}) and equate the coefficients of
$r^n$ on both sides of (\ref{eq:11}). Straightforward computations
show that
\begin{eqnarray}
\label{eq:iso7}
&& a_1 =a_3 = 0, ~\textrm{$a_2>0$ is arbitrary} \nonumber \\
&& a_4 = -\frac{a_2}{14} \nonumber \\ &&
h^*(r) = a_2\left[r^2 - \frac{r^4}{14}+\ldots\right]
\end{eqnarray}
where $a_2>0$ since $h^*$ is non-negative from Proposition~\ref{prop:1}.

Next, we show that the formal expansion (\ref{eq:iso6}) involves no odd powers of $r$. Direct computations show that $a_1 = a_3 =0$, as stated in (\ref{eq:iso7}). We proceed by induction. Suppose that $a_{2n+1}=0$ for $n=0\ldots p$. We show that $a_{2p+3}=0$ too.
Consider the left-hand side of the ordinary differential equation (\ref{eq:11}) i.e.
$$ \frac{d^2 h^*}{ dr^2} + \frac{2}{r}\frac{d h^*}{d r} - \frac{ 6 h^*}{r^2} =
\sum_{p=0}^{\infty} r^{n-2} a_{n} \left[n^2 + n - 6\right] $$ so
that the coefficient of $r^{2p+1}$ is $(4p+2)(p+3)a_{2p+3}$. We
compute the coefficient of $r^{2p+1}$ on the right-hand side of
(\ref{eq:11}). One can directly show that
$$
-h^*+ h{^{*}}^3 + \frac{3 h_+}{t}\left( h{^{*}}^3- h{^{*}}^2\right) = \sum_{n=1}^{\infty} b_n r^n
$$
where
\begin{equation}
\label{eq:isonew} b_{2p+1} = - a_{2p+1} + \left(1 +
\frac{3h_+}{t}\right) \left[ 3(a_1^2 a_{2p-1} + a_2^2 a_{2p-3} +
\ldots + a_p^2 a_1) + a^3_{\frac{2p+1}{3}}\right] -
\frac{6h_+}{t}\left( a_1 a_{2p} + a_2 a_{2p-1} + \ldots + a_{p}
a_{p+1}\right)
\end{equation}
where the term involving $a_{\frac{2p+1}{3}}$ comes into play if
$\frac{2p+1}{3}$ is a positive integer. One can check
(\ref{eq:isonew}) by noting that the coefficient of $r^{2p+1}$ in
the series $h{^{*}}^2$ is $\sum_{n=1}^{2p} 2 a_n a_{2p+1 -n}$ so
that both $\left\{n, 2p+1-n \right\}\leq 2p+1$ and one of
$\left\{n, 2p+1-n\right\}$ is odd. Similarly, we note that the
coefficient of $r^{2p+1}$ in the series $h{^{*}}^3$ is
$a^3_{\frac{2p+1}{3}} + \sum_{n=1}^{p} 3 a_{n}^2 a_{2p+1 - 2n}$,
$\left\{n, 2p+1-2n\right\} < 2p+ 1$, $\frac{2p+1}{3}$ and
$2p+1-2n$ are necessarily odd. However, from the hypothesis,
$a_{2n+1}=0$ for $n=0\ldots p$. Therefore, $b_{2p+1}=0$ in
(\ref{eq:isonew}) and since
$$ b_{2p+1} = (4p+2)(p+3)a_{2p+3},$$
we deduce that $a_{2p+3}=0$ as required.

The bounds (\ref{eq:bounds}) have been established in
\cite{mkaddem&gartland1,mkaddem&gartland2} and are valid in the
$R\to\infty$ limit. The inequalities (\ref{eq:bounds2}) follow
from (\ref{eq:bounds}) and the limit
$$ a_2 = \lim_{r \to 0} \frac{h^*(r)}{r^2}.$$  $\Box$

\begin{prop} \label{prop:biaxial}
Consider the radial-hedgehog solution
$$\Q^*(\rvec) = \sqrt{\frac{3}{2}}h^*(r)\left(\frac{\rvec\otimes \rvec}{r^2}  -\frac{1}{3}\I\right)$$
where $h^*$ is a global minimizer of $I$ in (\ref{eq:10}) in the
admissible space $\Acal_h$. Then $\Q^*$ is not the global
minimizer of $\tilde{\mathcal{I}}_{LG}$ in (\ref{eq:nondim4}) in
the admissible space \begin{equation} \label{eq:AQ} \Acal_\Q =
\left\{\Q \in W^{1,2}\left(B(0,R);\bar{S}\right);~ \Q =
\sqrt{\frac{3}{2}}\left(\frac{\rvec}{r}\otimes \frac{\rvec}{r} -
\frac{\mathbf{I}}{3}\right)~on~\partial B(0,R)\right\}
\end{equation} (which is the
admissible space defined in (\ref{eq:5}) in terms of the
dimensionless variables in (\ref{eq:nondim}) and
(\ref{eq:nondimnew})), for sufficiently large values of the domain
size $R$ and the reduced temperature $t$. In particular, the
biaxial state
\begin{equation}
\label{eq:biaxialperturbation} \hat{\Q}(\rvec) = \begin{cases}
\Q^*(\rvec) + \frac{1}{(r^2 + 12)^2}\left(1 -
\frac{r}{\sigma}\right)\left(\zvec\otimes\zvec -
\frac{\mathbf{I}}{3}\right) ~ 0\leq r\leq \sigma \\
\Q^*(\rvec)~\sigma \leq r\leq R ,
\end{cases}
\end{equation}
where $\zvec=(0,0,1)$ is the unit-vector in the $z$-direction, has lower free energy than the radial-hedgehog solution for
\begin{equation}
\label{eq:sigma}
\sigma = 10 \quad t=200.
\end{equation}
\end{prop}

\textit{Proof:}
Consider a general biaxial perturbation (\ref{eq:biaxialperturbation})
$$
\hat{\Q}(\rvec) = \begin{cases} \Q^*(\rvec) +
p(r)\left(\zvec\otimes\zvec - \frac{\mathbf{I}}{3}\right) ~ 0\leq
r\leq \sigma \\ \Q^*(\rvec)~\sigma \leq r\leq R ,
\end{cases}
$$ where $p:[0,R]\to \mathbb{R}$ is non-zero for $0\leq r\leq \sigma$ and $p(r)=0$ for all $\sigma\leq r\leq R$, $\Qvec^*$ is the radial-hedgehog solution in (\ref{eq:13}), $\rvec = \left(x,y,z\right)$ is the position vector, $\zvec=(0,0,1)$ is the unit-vector
in the $z$-direction and $\mathbf{I}$ is the $3\times 3$ identity matrix. In particular, the perturbation $\bar{\Q}$ is localised in a ball of radius $\sigma$ around the origin or equivalently, is localised around the isotropic core of the radial-hedgehog solution and the radius $\sigma$ will be determined as part of the problem.

Let $(r,\theta,\phi)$ with $r\in[0,R], \theta\in\left[0,\pi\right], \phi\in[0,2\pi)$ denote a spherical coordinate system centered at the origin.
Straightforward computations show that
\begin{eqnarray}
\label{eq:b4}
&& |\grad \hat{\Q}|^2 = |\grad \Q^*|^2 + \frac{2}{3}\left(\frac{dp}{dr}\right)^2 + \sqrt{6}\frac{d h^*}{dr}\frac{d p}{dr} \left(\cos^2\theta - \frac{1}{3}\right) \nonumber \\
&& \textrm{tr}\hat{\Q}^2 = \textrm{tr}\Q{^*}^2 + \frac{2}{3}p^2(r) + \sqrt{6} h(r) p(r)\left(\cos^2\theta - \frac{1}{3}\right) \nonumber \\
&& \textrm{tr}\hat{\Q}^3 =
\textrm{tr}\Q{^*}^3 + \frac{2}{9}p^3(r) + \left(\sqrt{\frac{2}{3}} + \frac{1}{\sqrt{6}}\right)h^*(r)p^2(r)\left(\cos^2\theta - \frac{1}{3}\right) + \frac{3}{2}h{^*}^2(r) p(r)\left(\cos^2\theta - \frac{1}{3}\right) \nonumber \\
&& \left(\textrm{tr}\hat{\Q}^2\right)^2 =
\left(\textrm{tr}\Q{^*}^2 \right)^2 + \frac{4}{9}p^4(r) + 6
\left(h^*(r)\right)^2 p^2(r)\left(\cos^2\theta -
\frac{1}{3}\right)^2 + \nonumber \\ && +
\frac{4}{3}\left(h^*(r)\right)^2 p^2(r) + 2\sqrt{6}
h^*(r)p(r)\left(\cos^2\theta - \frac{1}{3}\right)\left[h{^*}^2 +
\frac{2}{3}p^2(r)\right].
\end{eqnarray}

Noting that $$\int_{0}^{\pi}\left(\cos^2\theta - \frac{1}{3}\right)\sin\theta d\theta = 0$$ and
$$\int_{0}^{\pi}\left(\cos^2\theta - \frac{1}{3}\right)^2\sin\theta d\theta = \frac{8}{45},$$ we obtain the following -
\begin{eqnarray}
&& \frac{1}{4\pi}\left[\tilde{\mathcal{I}}_{LG}[\hat{\Q}] - \tilde{\mathcal{I}}_{LG}[\Q^*]\right] = \nonumber \\
&& = \int_{0}^{\sigma} \frac{r^2}{3}\left(\frac{d p}{dr}\right)^2 - \frac{r^2}{3} p^2(r) - 2\sqrt{6}\frac{h_+}{9t}r^2 p^3(r) + \frac{r^2 h_+^2}{2t}\left\{\frac{4}{9}p^4(r)+\frac{28}{15}h{^*}^2 p^2(r)\right\}~dr
\label{eq:b8}
\end{eqnarray} where $h_+$ has been defined in (\ref{eq:h+}).
Recalling the bounds (\ref{eq:bounds}), we have that
\begin{eqnarray}
&& \frac{1}{4\pi}\left[\tilde{\mathcal{I}}_{LG}[\hat{\Q}] - \tilde{\mathcal{I}}_{LG}[\Q^*]\right] < \nonumber \\
&& < \int_{0}^{\sigma} \frac{r^2}{3}\left(\frac{d p}{dr}\right)^2 - \frac{r^2}{3} p^2(r) - 2\sqrt{6}\frac{h_+}{9t}r^2 p^3(r) + \frac{r^2 h_+^2}{2t}\left\{\frac{4}{9}p^4(r)+\frac{28}{15}\left(\frac{r^2}{r^2 + t \lambda_t^2}\right)^2 p^2(r)\right\}~dr
\label{eq:b9}
\end{eqnarray} where $\lambda_t^2 = \frac{24}{9 + 8t + 3\sqrt{9+8t}}$.
Let \begin{equation}\label{eq:p} p(r) =  \frac{1}{(r^2 + 12)^2}\left(1 - \frac{r}{\sigma}\right).\end{equation}
One can then directly substitute (\ref{eq:p}), $\sigma=10$ and $t=200$ into (\ref{eq:h+}) and (\ref{eq:b9}) to find that the associated free energy difference
$$\frac{1}{4\pi}\left[\tilde{\mathcal{I}}_{LG}[\hat{\Q}] - \tilde{\mathcal{I}}_{LG}[\Q^*]\right] < 0$$
i.e. we have found a biaxial perturbation localised in a ball $B(0,\sigma)$, that has lower free energy than the radial-hedgehog solution for $t=200$. Therefore, the radial-hedgehog solution cannot be a global Landau-de Gennes minimizer in this regime. Proposition~\ref{prop:iso2} now follows. $\Box$

The instability of the radial-hedgehog solution with respect to
biaxial perturbations has been theoretically demonstrated in
\cite{mkaddem&gartland2}, in the limit $t\to\infty$ and
$R\to\infty$. The authors in \cite{mkaddem&gartland2} consider the
second variation of the Landau-de Gennes energy functional and
treat the instability condition as a Schrodinger eigenvalue
problem, which has to be solved numerically. We construct an
explicit biaxial perturbation, localized near the isotropic core
of the radial-hedgehog solution and show that this biaxial
perturbation has lower free energy than the radial-hedgehog
solution, for sufficiently low temperatures. The biaxial
perturbation is energetically preferable only when localized in a
ball $B(0,\sigma)$ centered at the origin and one can check that
$\left[\tilde{\mathcal{I}}_{LG}[\hat{\Q}] -
\tilde{\mathcal{I}}_{LG}[\Q^*]\right]>$ if $\sigma$ is too small
or too large i.e. $\sigma$ needs to be large enough for the
biaxiality to manifest itself and yet be small enough so as not to
perturb the far-field properties. We demonstrate instability for
$t=200$ and $\sigma=10$ and the parameter regimes can be
investigated more systematically. In particular, our approach in
Proposition~\ref{prop:iso2} gives insight into how to quantify the
instability regime analytically.

\section{The limit $t\to\infty$}
\label{sec:lowtemp}

Consider the ordinary differential equation in (\ref{eq:11})
$$\frac{d^2 h}{dr^2} + \frac{2}{r}\frac{dh}{dr} - \frac{6h}{r^2} = -h + h^3 + \frac{3 h_+}{t}\left(h^3 - h^2\right)$$
in the limit $t\to\infty$. In the limit $t\to\infty$, $$ \frac{h_+}{t}\leq \frac{\beta}{\sqrt{t}}$$ for some $\beta>0$ independent of $t$ and hence for any non-negative solution $h$, we have
$$\left|\frac{3 h_+}{t}\left(h^3 - h^2\right)\right| << h - h^3,$$
 since $0\leq h(r)\leq 1$. Recall that the upper bound in (\ref{eq:15}) applies to all solutions of (\ref{eq:11}) and not just the radial-hedgehog solution. In the limit $t \to \infty$, the ordinary differential equations (\ref{eq:11}) approximately reduces to
 \begin{equation}
 \label{eq:lowtemp}
 \frac{d^2 h}{dr^2} + \frac{2}{r}\frac{dh}{dr} - \frac{6h}{r^2} \approx -h + h^3
 \end{equation}
 although the influence of the perturbation term $\frac{3 h_+}{t}\left(h^3 - h^2\right)$ needs to be carefully quantified. The limiting problem (\ref{eq:lowtemp}) has a very similar structure to the governing ordinary differential equation for vortex solutions in the Ginzburg-Landau theory of superconductivity \cite{bbh1}. Vortex-solutions have been widely studied within the Ginzburg-Landau framework \cite{farina,herve&herve}. They have the special structure
 $$ w(\xvec) = u(|\xvec|) g\left(\frac{\xvec}{|\xvec|}\right) \quad \xvec\in\mathbb{R}^N$$
 where $u$ is a solution of the following ordinary differential equation in $\mathbb{R}^N$
 \begin{eqnarray}\label{eq:lowtemp2}
&& \frac{d^2 u}{d|\xvec|^2} + \frac{N-1}{|\xvec|} - \frac{\lambda_K}{|\xvec|^2}u =  -u + u^3 \nonumber \\
 && u(0) = 0
 \end{eqnarray}
 and $\lambda_K$ is a characteristic constant. In what follows, we adapt
 Ginzburg-Landau techniques for (\ref{eq:lowtemp2}) to the ordinary
 differential equation (\ref{eq:11}) in the limit $t\to\infty$ to establish
 uniqueness and global monotonicity of $h^*$ in (\ref{eq:13}). As will be demonstrated
 in Section~\ref{sec:ns}, there are important technical differences between (\ref{eq:11}) and the Ginzburg-Landau formulation (\ref{eq:lowtemp2}) and in general, Ginzburg-Landau results do not readily transfer to the Landau-de Gennes framework. In this sense, one could also refer to the $t\to\infty$ limit as the \emph{Ginzburg-Landau limit}.

\begin{lem} \cite{bbh2}
\label{lem:gradient}
For all $t>1$ and any solution $\Q$ of the Euler-Lagrange equations (\ref{eq:7}), we have the following global upper bound for the gradient -
\begin{equation}
\label{eq:gradient}
\| \grad \Q \|_{L^{\infty}(B(0,R))} \leq C
\end{equation}
where $C>0$ is independent of $t$. For the radial-hedgehog
solution $\Q^*(\rvec) = \sqrt{\frac{3}{2}}
h^*(r)\left(\hat{\rvec}\otimes\hat{\rvec} -
\frac{1}{3}\mathbf{I}\right)$, this implies that for $t>1$, we
have the following inequality
\begin{equation}
\label{eq:gradient2}
\left|\grad \Q^*\right|^2 = \left(\frac{ d h^*}{ d r}\right)^2 + \frac{ 3 h{^*}^2}{r^2} \leq C^2 \quad \forall r \in [0,R]
\end{equation} where $C$ is again independent of $t$.
\end{lem}

\textit{Proof:} The proof of Lemma~\ref{lem:gradient} can be found in \cite{bbh2} where the authors show that a solution $u$ of the elliptic system
$$-\Delta u = f \quad ~on ~\Omega\subset \mathbb{R}^n $$
satisfies
$$|\grad u (\rvec)|^2 \leq C \| f \|_{L^{\infty}(\Omega)}\| u \|_{L^{\infty}(\Omega)} \quad \forall\rvec\in\Omega.$$
In our case, we apply this result to the system (\ref{eq:7}), noting that $\Q^*$ is a solution of (\ref{eq:7}),
$$ f =  - \Q_{ij}  - \frac{ 3\sqrt{6} h_+}{t}\left(\Q_{ik}\Q_{kj}-\frac{\delta_{ij}}{3}\textrm{tr}(\Q^2)\right)
+\frac{2 h_+^2}{t}\Q_{ij}\left(\textrm{tr}\Q^2\right) $$ for each $i,j=1\ldots 3$, $\| \Q^*\|_{L^{\infty}(\Omega)} \leq 1$ from the bounds in (\ref{eq:15})
and
$\frac{h_+}{t}\leq \frac{9}{4}$ and $\frac{h_+^2}{t} \leq \frac{33}{4}$ for $t>1$. $\Box$

\begin{lem}
\label{lem:far-field}
In the limit $t\to \infty$, we have that
\begin{equation}
\label{eq:far}
\lim_{r\to \infty} r^2 \frac{ d h^*}{dr} = 0.
\end{equation}
\end{lem}
\textit{Proof:}
In the limit $t\to \infty$, we solve the ordinary equation (\ref{eq:11})
on an unbounded domain i.e. the boundary conditions (\ref{eq:12}) become
\begin{equation}
\label{eq:bc}
h(0) = 0 \quad h(r)\to 1~as~ r\to \infty.
\end{equation}
From the bounds (\ref{eq:bounds}), we deduce that for $r$ sufficiently large,
\begin{equation}
\label{eq:bc2} h^*(r) = 1+ \sigma (r) \quad ~where~ -
\frac{\alpha}{r^2} \leq \sigma(r) \leq -\frac{\beta}{r^2}
\end{equation}
for positive constants $\alpha, \beta$ independent of $t$. These
bounds imply that $h^*(r) \to 1$ uniformly as $r\to\infty$ and
from Proposition~\ref{prop:farfield} in the next section, this
implies that $$\frac{dh^*}{dr} = \frac{d\sigma}{dr}>0$$ for $r$
sufficiently large.

We use (\ref{eq:11}) to obtain an ordinary differential equation
for $\delta = \frac{dh^*}{dr}$ as shown below :
\begin{equation}
\label{eq:bc3} r^2 \frac{d^2 \delta}{d r^2} =
-4r\frac{d\delta}{dr} + 4\delta + 2r\left(h{^*}^3 - h^*\right) +
r^2\left(3 h{^*}^2 -1 \right)\delta
\end{equation}
where $\delta \to 0$ as $r\to\infty$. We can then use differential
inequalities as in \cite{mkaddem&gartland2} to deduce that
\begin{equation}
\label{eq:bc4} \delta = \frac{dh^*}{dr} \leq \frac{\gamma_1}{r^3}
\end{equation}
where $\gamma > 2 \beta>0$ is a positive constant and $\beta$ has
been in defined in (\ref{eq:bc2}). Since $\beta$ is independent of
$t$ and $r$, $\gamma$ can be chosen to be a positive constant
independent of $t$ and (\ref{eq:bc4}) implies that
\begin{equation}
\label{eq:bc5} \lim_{r\to\infty} r^2\frac{dh^*}{dr} = 0.
\end{equation} Lemma~\ref{lem:far-field} now follows.$\Box$

\begin{prop}
\label{prop:uniqueness} The ordinary differential equation
\begin{equation}
\label{eq:new}
\frac{d^2 h}{dr^2} + \frac{2}{r}\frac{dh}{dr} - \frac{6h}{r^2} = -h + h^3 + \frac{3 h_+}{t}\left(h^3 - h^2\right)
\end{equation}
subject to the boundary conditions
\begin{equation}
\label{eq:new2}
h(0) = 0 \quad h(r)\to 1 \quad r\to\infty
\end{equation}
has a unique non-negative solution in the limit $t\to \infty$.
\end{prop}

\textit{Proof:} Let $h_1$ and $h_2$ be two different non-negative solutions of (\ref{eq:new}) subject to the boundary conditions (\ref{eq:new2}) i.e.
\begin{eqnarray}
&& \frac{h_{1}^{''}}{h_1} + \frac{2}{r}\frac{h_1^{'}}{h_1} - \frac{6}{r^2} + \left(1 - h_1^2\right) + \frac{ 3 h_+}{t} h_1 \left(1 - h_1^2\right) = 0 \nonumber\\ && \frac{h_{2}^{''}}{h_2} + \frac{2}{r}\frac{h_2^{'}}{h_2} - \frac{6}{r^2} + \left(1 - h_2^2\right) + \frac{ 3 h_+}{t} h_2 \left(1 - h_2^2\right) = 0
\end{eqnarray} where $h_1^{'} = \frac{ d h_1}{ dr}$, $h_1^{''} = \frac{d^2 h_1}{d r^2}$ etc.
We subtract the two equations to get
\begin{equation}
\label{eq:u1}
 - \frac{h_{1}^{''}}{h_1} + \frac{h_{2}^{''}}{h_2} - \frac{2}{r}\left(\frac{h_1^{'}}{h_1} - \frac{h_2^{'}}{h_2}\right) =
 \left(1 + \frac{3 h_+}{t}\right)\left( h_2^2 - h_1^2\right) + \frac{ 3 h_+}{t}\left( h_1 - h_2 \right).
 \end{equation}
 Following the methods in \cite{bbh1}, we multiply both sides of (\ref{eq:u1}) by
 $r^2\left(h_1^2 - h_2^2 \right)$ and integrate from $r=0$ to $r=R$ to find
 \begin{eqnarray}
 \label{eq:u2}
 && \int_{0}^{R} r^2 \left( \frac{h_1}{h_2} h_2^{'} - h_{1}^{'} \right)^2~ dr + \int_{0}^{R} r^2 \left( \frac{h_2}{h_1} h_1^{'} - h_{2}^{'} \right)^2~ dr
 + \nonumber \\ && + \int_{0}^{R}\left( 1 + \frac{ 3 h_+}{t}\right) r^2 \left( h_1^2 - h_2^2 \right)^2~dr - \frac{ 3 h_+}{t}\int_{0}^{R}\left(h_1 - h_2\right)^2 r^2 \left(h_1 + h_2\right)~dr = \nonumber
  \\ && = -r^2 h_2^{'}\left(\frac{h_1^2}{h_2} - h_2 \right)\Large{\mid}_{0}^{R}
  - r^2 h_{1}^{'}\left(\frac{h_2^2}{h_1} - h_1 \right)\Large{\mid}_{0}^{R}.
 \end{eqnarray}
 Taking the limit $R \to \infty$ and using (\ref{eq:far}), we have that
 \begin{eqnarray}
 \label{eq:u3}
&&  \lim_{R \to \infty} \int_{0}^{R} r^2 \left( \frac{h_1}{h_2} h_2^{'} - h_{1}^{'} \right)^2~ dr + \int_{0}^{R} r^2 \left( \frac{h_2}{h_1} h_1^{'} - h_{2}^{'} \right)^2~ dr + \nonumber \\ && +
  \int_{0}^{R}\left( 1 + \frac{ 3 h_+}{t}\right) r^2 \left( h_1^2 - h_2^2 \right)^2~dr - \frac{ 3 h_+}{t}\int_{0}^{R}\left(h_1 - h_2\right)^2 r^2 \left(h_1 + h_2\right)~dr = 0.
 \end{eqnarray}
 From (\ref{eq:u3}), we deduce that
\begin{equation}
\label{eq:unew}\int_{0}^{R}\left( 1 + \frac{ 3 h_+}{t}\right) r^2 \left( h_1^2 - h_2^2 \right)^2~dr - \frac{ 3 h_+}{t}\int_{0}^{R}\left(h_1 - h_2\right)^2 r^2 \left(h_1 + h_2\right)~dr \to 0 \quad R\to\infty.
\end{equation}

 We first make the elementary observation that $\exists R_1 \in [0,R]$ such that
 \begin{equation}
 \label{eq:u4}
 h_1 (r ) + h_2 (r) > 1 \quad \forall r > R_1.
 \end{equation}
 The inequalities (\ref{eq:bounds}) are true for any solution of (\ref{eq:new}) subject to the boundary conditions (\ref{eq:new2}) \cite{mkaddem&gartland1,mkaddem&gartland2}. Therefore, we have for $t>1$,
 \begin{equation}
 \label{eq:u5}
 \frac{r^2}{r^2 + 14} \leq h_1(r), h_2(r) \leq \frac{r^2}{r^2 + 6/17}
 \end{equation}
 and
 \begin{equation}
 \label{eq:R1}
 \frac{6}{17}\leq R_1^2 \leq 14
 \end{equation} i.e. $R_1$ can be bounded independently of $t$ for $t>1$.

 We partition the integral contribution in (\ref{eq:u3}) into two sub-intervals $[0,R_1]$ and $[R_1,R]$ respectively.
 \begin{eqnarray}
 \label{eq:u6}
 && \int_{0}^{R}\left( 1 + \frac{ 3 h_+}{t}\right) r^2 \left( h_1^2 - h_2^2 \right)^2~dr -
 \frac{ 3 h_+}{t}\int_{0}^{R}\left(h_1 - h_2\right)^2 r^2 \left(h_1 + h_2\right)~dr = \nonumber
 \\ && =\int_{0}^{R_1}\left( 1 + \frac{ 3 h_+}{t}\right) r^2 \left( h_1^2 - h_2^2 \right)^2 -
 \frac{3h_+}{t}r^2\left(h_1 - h_2\right)^2  \left(h_1 + h_2\right)~dr + \nonumber \\ && + \int_{R_1}^{R} r^2 \left( h_1^2 - h_2^2 \right)^2 + r^2\frac{3 h_+}{t}\left(h_1 - h_2\right)^2  \left(h_1 + h_2\right)\left[ h_1 + h_2 - 1\right]~dr
 \end{eqnarray}
 and note from (\ref{eq:u4}) that
 $$  \int_{R_1}^{R} r^2 \left( h_1^2 - h_2^2 \right)^2
 + r^2\frac{3 h_+}{t}\left(h_1 - h_2\right)^2  \left(h_1 + h_2\right)\left[ h_1 + h_2 - 1\right]~dr
 \geq \int_{R_1}^{R}r^2 \left( h_1^2 - h_2^2 \right)^2~dr.$$

\textit{Claim: For $t$ sufficiently large,}
\begin{equation}
\label{eq:u7} \int_{0}^{R_1}\left( 1 + \frac{ 3 h_+}{t}\right) r^2
\left( h_1^2 - h_2^2 \right)^2 - \frac{3h_+}{t}r^2\left(h_1 -
h_2\right)^2 \left(h_1 + h_2\right)~dr >
\frac{1}{2}\int_{0}^{R_1}r^2 \left( h_1^2 - h_2^2 \right)^2 ~dr.
\end{equation}

Recalling that $R_1$ can be bounded independently of $t$, we note that
$$
\frac{3h_+}{t} \int_{0}^{R_1} r^2\left(h_1 - h_2\right)^2  \left(h_1 + h_2\right)~dr \leq \gamma_1 \frac{h_+}{t}R_1^3 \leq \frac{\gamma_2}{\sqrt{t}}
$$
where $\gamma_1$ and $\gamma_2$ are positive constants independent of $t$. Therefore, the claim in (\ref{eq:u7}) is equivalent to
\begin{equation}
\label{eq:u8}
\sqrt{t} \geq \frac{\gamma_3}{\int_{0}^{R_1} r^2 \left(h_1^2  -  h_2^2 \right)^2~dr}
\end{equation}
for a positive constant $\gamma_3$ independent of $t$.

We note that $$\int_{0}^{R_1} r^2 \left(h_1^2  -  h_2^2 \right)^2~dr \leq \frac{R_1^3}{3}$$
so that as $t\to\infty$, we have two possibilities - (a) $\int_{0}^{R_1} r^2 \left(h_1^2  -  h_2^2 \right)^2~dr = O(1)$ as $t\to\infty$ and (b) $\int_{0}^{R_1} r^2 \left(h_1^2  -  h_2^2 \right)^2~dr  = o(1)$ as $t\to\infty$. In case (a), the condition (\ref{eq:u8}) is clearly satisfied for $t$ sufficiently large and the claim (\ref{eq:u7}) follows.

For case (b), we have
\begin{equation}
\label{eq:u10}\int_{0}^{R_1} r^2 \left(h_1^2  -  h_2^2 \right)^2~dr \to 0 ~as~t\to\infty.
\end{equation} From Lemma~\ref{lem:gradient} and the global bound (\ref{eq:15}), we obtain
\begin{equation}
\label{eq:u9}
\left|\grad\left(h_1^2 - h_2^2\right)\right|\leq D
\end{equation}
where $D$ is a positive constant independent of $t$.
Consider $r_0 \in [0,R_1]$ and let
$$ \left|\left(h_1^2 - h_2^2\right)(r_0)\right| = \alpha_0 >0.$$
Then from (\ref{eq:u9}), we have that
$$\left|\left(h_1^2 - h_2^2\right)(r)\right| \geq \frac{\alpha_0}{2} \quad r\in\left[r_0 - \frac{\alpha_0}{2D}, r_0 + \frac{\alpha_0}{2D}\right] $$
and therefore
$$ \int_{0}^{R_1} r^2 \left(h_1^2  -  h_2^2 \right)^2~dr \geq \int_{r_0 - \frac{\alpha_0}{2D}}^{r_0 + \frac{\alpha_0}{2D}}\frac{\alpha_0^2}{4} r^2~dr \geq \gamma_4 \alpha_0^5$$
where $\gamma_4$ is a positive constant independent of $t$.
Combining the above with (\ref{eq:u10}), we have that $\alpha_0 \to 0$ as $t\to \infty$ and hence
\begin{equation}
\label{eq:u11}
h_1(r) = h_2(r) \quad r\in\left[0,R_1\right]
\end{equation}
since the choice of $r_0$ is arbitrary and we are interested in non-negative solutions.

From (\ref{eq:unew}) and (\ref{eq:u4}), we have that
\begin{eqnarray}
\label{eq:u12}
&& \int_{0}^{R}\left( 1 + \frac{ 3 h_+}{t}\right) r^2 \left( h_1^2 - h_2^2 \right)^2~dr - \frac{ 3 h_+}{t}\int_{0}^{R}\left(h_1 - h_2\right)^2 r^2 \left(h_1 + h_2\right)~dr \geq \nonumber \\ && \geq \int_{0}^{R_1} \left( 1 + \frac{ 3 h_+}{t}\right) r^2 \left( h_1^2 - h_2^2 \right)^2~dr - \frac{ 3 h_+}{t}\int_{0}^{R_1}\left(h_1 - h_2\right)^2 r^2 \left(h_1 + h_2\right)~dr + \int_{R_1}^{R}r^2 \left( h_1^2 - h_2^2 \right)^2~dr.
\end{eqnarray}
For case (a),(\ref{eq:u7}) holds and (\ref{eq:u12}) can be written as
\begin{eqnarray}
&& \frac{1}{2}\int_{0}^{R_1}r^2 \left( h_1^2 - h_2^2 \right)^2 ~dr + \int_{R_1}^{R}r^2 \left( h_1^2 - h_2^2 \right)^2~dr \leq \nonumber \\ && \leq \int_{0}^{R}\left( 1 + \frac{ 3 h_+}{t}\right) r^2 \left( h_1^2 - h_2^2 \right)^2~dr - \frac{ 3 h_+}{t}\int_{0}^{R}\left(h_1 - h_2\right)^2 r^2 \left(h_1 + h_2\right)~dr \to 0 \quad t\to\infty
\end{eqnarray}
from which we deduce that
\begin{equation}
\left( h_1^2 - h_2^2 \right)^2 = 0 \quad r\in \left[0, R\right]
\end{equation}
or equivalently
\begin{equation}
h_1(r) = h_2(r) \quad r\in \left[0, R\right]
\end{equation}

For case (b), we have established in (\ref{eq:u11}) that
$h_1(r) = h_2(r) \quad r\in\left[0,R_1\right]$ and hence
$$\int_{0}^{R_1}
\left( 1 + \frac{ 3 h_+}{t}\right) r^2 \left( h_1^2 - h_2^2
\right)^2~dr - \frac{ 3 h_+}{t}\int_{0}^{R_1}\left(h_1 -
h_2\right)^2 r^2 \left(h_1 + h_2\right)~dr \to 0 \quad
t\to\infty.$$ From (\ref{eq:unew}), we deduce that
$$ \int_{R_1}^{R}r^2 \left( h_1^2 - h_2^2 \right)^2~dr \to 0 \quad t \to\infty$$
and hence,
$$ h_1 (r) = h_2(r) \quad r\in \left[R_1, R \right.]$$
Combining the above with (\ref{eq:u11}), we have that
\begin{equation}
h_1(r) = h_2(r) \quad r\in \left[0, R\right]
\end{equation} in case (b) too. Proposition~\ref{prop:uniqueness} now follows. $\Box$

We, next, illustrate the applicability of shooting arguments to
the ordinary differential equation (\ref{eq:new}) in the limit
$t\to\infty$, so that the corresponding boundary conditions are
(\ref{eq:new2}) \cite{tang, farina}.

From Proposition~\ref{prop:iso2}, we have that  for any solution
$h$ of (\ref{eq:new}) subject to the boundary condition
$$ h( 0) =0,$$ $\exists$ a constant $a_2$ such that
\begin{equation}
\label{eq:shoot}
h(r) \sim a_2 r^2 \quad r\to 0.
\end{equation} Given $a_2$, we denote the corresponding solution by $h(a_2,r)$.
We are interested in non-negative solutions and hence, we take $a_2 >0$. By analogy with \cite{tang,farina}, we call $a_2$ the \emph{shooting parameter}. We consider three different classes of solutions
\begin{itemize}
\item $\mathcal{P} = \left\{ a_2 >0 ;~ \exists z \in \left(0, R_a\right) ~such~that~ \frac{d h(a_2,z)}{dr} = 0\right\}$
\item $\mathcal{Q} = \left\{a_2 >0 ;~ \frac{d h(a_2,z)}{dr} > 0~and~ h(a_2,r)\leq 1~for~all~ r> 0 \right\}$
\item $\mathcal{R} = \left\{a_2 >0 ;~ \frac{d h(a_2,z)}{dr} > 0 \forall r\in\left(0, R_a\right)~and~ \max_{r\in\left(0, R_a\right)} h\left(a_2, r \right)>1 \right\}$
\end{itemize} where $R_a$ is the maximal interval of existence of the solution $h(a_2,r)$.
Clearly
$$\mathcal{P}\cap \mathcal{Q} = \mathcal{P}\cap \mathcal{R} = \mathcal{Q}\cap \mathcal{R}= \phi $$
and
$$ \mathcal{P}\cup\mathcal{Q}\cup\mathcal{R} = (0,\infty).$$

Our aim is to show that $\mathcal{P}$ and $\mathcal{R}$ are non-empty and open. Then, $\mathcal{Q}$ is also non-empty. We have a unique solution of the ordinary differential equation (\ref{eq:new}) subject to the boundary conditions (\ref{eq:new2}) in the limit $t\to\infty$. Therefore, if we can show that $a_2 \in \mathcal{Q}$ implies that the corresponding $h(a_2, r)$ is a solution of (\ref{eq:new}) and (\ref{eq:new2}), then we have that $h^* \in \mathcal{Q}$ in the limit $t\to\infty$ and hence, $\frac{d h^*}{ dr}>0$ for all $r>0$ i.e. we have global monotonicity in the limit $t\to\infty$.

It is evident that a solution of (\ref{eq:new}) subject to the
boundary conditions (\ref{eq:new2}) cannot belong to $\mathcal{R}$
owing to the global bounds (\ref{eq:15}). It remains to rule out
the possibility $a_2\in\mathcal{P}$. We start with an elementary
lemma.

\begin{lem}
\label{lem:Q}
If $a_2\in\mathcal{Q}$, then $h(a_2,r)$ is a solution of (\ref{eq:new}) subject to the boundary conditions (\ref{eq:new2}).
\end{lem}
\textit{Proof:} The proof closely follows the methods in \cite{tang}. Since $h(a_2,r)$ is monotonically increasing (from the definition of $\mathcal{Q}$) and is bounded above by $1$, $b=\lim_{r\to\infty} h(a_2, r)$ exists and $b\in\left(0,1\right]$. Hence, to finish the proof, we need to show that $b=1$. In fact, if $b<1$, then as $r\to\infty$, (\ref{eq:11}) can be written as
$$ \frac{d}{dr}\left(r^2 \frac{dh}{dr}\right) = 6b + r^2\left( b^3 - b + \frac{ 3 h_+}{t}\left(b^3 - b^2\right)\right)$$
so that
$$\frac{dh}{dr}\sim \frac{6b}{r} + \frac{r}{3}\left(b^3 - b + \frac{3h_+}{t}\left(b^3 - b^2\right)\right)$$
contradicting the hypothesis that $\frac{dh}{dr}>0$ for all $r>0$. Therefore, $b=1$ and Lemma~\ref{lem:Q} follows. $\Box$

Next, we need to show that $\mathcal{Q}$ is non-empty. For this, we need
\begin{lem}
\label{lem:P} The set $\mathcal{P}$ is not empty; more precisely, there exists a positive constant $m$ such that $(0,m)\subset \mathcal{P}$.
\end{lem}
\textit{Proof:} The proof closely follows the methods in \cite{tang,farina}. For any $a_2>0$, set
\begin{equation}
\label{eq:w}
w(a_2,r) = \frac{h(a_2,r)}{a_2};
\end{equation}
then $w$ satisfies the following ordinary differential equation from (\ref{eq:new})
\begin{eqnarray}
\label{eq:w2}
&& \frac{d^2 w}{d r^2} + \frac{2}{r}\frac{dw}{dr} - 6\frac{w}{r^2} + w = a_2^2 w^3 + \frac{3 h_+}{t}\left(a_2^2 w^3 - a_2 w^2 \right) \nonumber \\
&& w\left(a_2, r \right) \sim r^2 \quad r\to 0.
\end{eqnarray}
Then as $a_2 \to 0$, $w(a_2,r)\to w(0,r)$ where $w(0,r)$ is the solution of
\begin{eqnarray}
\label{eq:w3}
&& \frac{d^2 w}{d r^2} + \frac{2}{r}\frac{dw}{dr} - 6\frac{w}{r^2} + w = 0 \nonumber \\
&& w(0,r)\sim r^2 \quad r\to 0
\end{eqnarray} and the general solution of this ordinary differential equation is
\begin{equation}
\label{eq:w4}
w(r) = \frac{C_1}{r^3}\left(-3\cos r - 3r\sin r + r^2\cos r \right) + C_2 r^2\left(-3\sin r + 3r \cos r + r^2 \sin r \right)
\end{equation}
for arbitrary constants $C_1$ and $C_2$. From (\ref{eq:w4}), we deduce that $w(a_2,r)$ has oscillatory behavior as $a_2 \to 0$ and hence, so does $h(a_2,r) = a_2 w(a_2,r)$. This completes the proof of the lemma. $\Box$

\begin{lem}
\label{lem:P2} The set $\mathcal{P}$ is open.
\end{lem}
\textit{Proof:} The proof of Lemma~\ref{lem:P2} closely follows
the methods in \cite{tang,farina} and we reproduce the proof for
completeness and to illustrate the technical differences.

For $a_2\in\mathcal{P}$, define
\begin{equation}\label{eq:m1}
z_0(a_2) = \inf\left\{r\in\left(0,R_a\right);~ \frac{
dh(a_2,r)}{dr} = 0\right\}
\end{equation} i.e. $z_0(a_2)$ is the smallest stationary point of $h(a_2,r)$. We can show that
\begin{equation}
\label{eq:m2}
\frac{ d^2 h \left(a_2, z_0(a_2)\right)}{ d r^2} < 0.
\end{equation}
The definition of $z_0(a_2)$ implies that
$$ \frac{ dh(a_2, z_0(a_2))}{dr} = 0~ and ~ \frac{ d^2 h \left(a_2, z_0(a_2)\right)}{ d r^2} \leq 0$$
since we are interested in non-negative solutions.

From the governing ordinary differential equation (\ref{eq:new}), we have that
\begin{equation}
\label{eq:m3}
\frac{d^2 h}{dr^2} + \frac{2}{r}\frac{dh}{dr}  = h\left( \frac{6}{r^2} + h^2 - 1 + \frac{3 h_+}{t}\left(h^2 - h\right)\right)\leq 0 ~at~ r=z_0(a_2)
\end{equation}
and note that
$$ \frac{d}{dr}\left[
\frac{6}{r^2} + h^2 - 1 + \frac{3 h_+}{t}\left(h^2 - h\right) \right] < 0 ~at ~r=z_0(a_2). $$

If $$\frac{ d^2 h \left(a_2, z_0(a_2)\right)}{ d r^2} = 0$$
then
$$ \frac{6}{r^2} + h^2 - 1 + \frac{3 h_+}{t}\left(h^2 - h\right) = 0 ~at ~r=z_0(a_2). $$
This implies that
$$\frac{6}{r^2} + h^2 - 1 + \frac{3 h_+}{t}\left(h^2 - h\right) > 0~on~ r\in\left[z_0(a_2) - \delta, z_0(a_2)\right)$$
for some $\delta>0$. On the other hand,
$$\frac{1}{r^2}\frac{d}{dr}\left[r^2 \frac{dh}{dr}\right] = h\left[\frac{6}{r^2} + h^2 - 1 + \frac{3 h_+}{t}\left(h^2 - h\right)\right]$$
from (\ref{eq:11}) so that
$$ \frac{1}{r^2}\frac{d}{dr}\left[r^2 \frac{dh}{dr}\right] >0 ~on~ r\in\left[z_0(a_2) - \delta, z_0(a_2)\right).$$
This in turn implies that
\begin{equation}
\label{eq:m4}
z_0^2(a_2)\frac{dh(a_2,z_0(a_2))}{dr} > \left(z_0(a_2) - \delta \right)^2 \frac{dh(a_2,z_0(a_2)-\delta)}{dr}>0
\end{equation}
(since $h_r>0$ for $r\in(0,z_0(a_2))$ from the definition (\ref{eq:m1})), contradicting the definition of $z_0(a_2)$. Hence, (\ref{eq:m2}) holds.

Finally, we note that for any $a_0\in\mathcal{P}$, by the Implicit Function Theorem and (\ref{eq:m2}), there exists a smooth function $y(a_2)$ defined in a neighbourhood of $a_0$ such that $y(a_0) = z_0(a_0)$ and $\frac{d h(a_2, y(a_2))}{dr} = 0$. Hence, $\mathcal{P}$ is open as required. $\Box$

\begin{lem} \label{lem:R}
The set $\mathcal{R}$ is non-empty and open.
\end{lem}
\textit{Proof:} The proof of Lemma~\ref{lem:R} closely follows the methods in \cite{farina}. We introduce the function
\begin{equation}
\label{eq:v}
v(r) = b h(a_2, br) ~ where~ b=a_2^{-1/3}.
\end{equation}
Then one can check that $v$ satisfies the following ordinary differential equation
\begin{equation}
\label{eq:v2}
\frac{d^2 v}{ dr^2} + \frac{2}{r}\frac{dv}{dr} - \frac{6v}{r^2} - \left(1 + \frac{3 h_+}{t}\right)v^3 + b^2 v + \frac{ 3 h_+}{t} b v^2 = 0
\end{equation}
with
\begin{equation}
\label{eq:v3}
v(r)\sim r^2 \quad r\to 0.
\end{equation}
If we let $b\to 0$, then the limiting problem is
\begin{eqnarray}
\label{eq:v4}
&& \frac{d^2 v}{ dr^2} + \frac{2}{r}\frac{dv}{dr} - \frac{6v}{r^2} - \left(1 + \frac{3 h_+}{t}\right)v^3  = 0 \nonumber \\
&& v(r)\sim r^2 \quad r\to 0.
\end{eqnarray}
From the hypothesis, we have that $v, \frac{dv}{dr}>0$ for $r>0$. We claim that there does not exist $l>0$ such that $\lim_{r\to\infty}v(r)=l$.
We prove the claim by contradiction. Assume $\exists l>0$ such that $\lim_{r\to\infty} v(r)=l$. Then (\ref{eq:v4}) implies that
$$ \frac{d}{dr}\left(r^2 \frac{dv}{dr}\right)\sim 6l + r^2\left(1+\frac{3h_+}{t}\right)l^3 \quad r\to\infty $$
so that
$$ \frac{ dv}{dr}\sim \frac{ 6 l}{r} + \frac{r}{3}\left(1+\frac{3h_+}{t}\right) l^3 \quad r\to\infty.$$
Therefore, $v(r)>>l$ for $r$ sufficiently large, which contradicts the hypothesis. The other possibility is $l=0$ but this contradicts the definition of $\mathcal{R}$ which requires that $\frac{dv}{dr}>0$ for all $r>0$. Therefore
\begin{equation}
\label{eq:v5}
v(r)\to \infty~\quad~as~r\to\infty.
\end{equation}
Consequently, $h(a_2,r)$ is large when $a_2$ is large enough and the set $\mathcal{R}$ is non-empty. By the continuous dependence of $h$ on $a_2$ and the definition of $\mathcal{R}$, we deduce that $\mathcal{R}$ is open. $\Box$
\begin{lem}
\label{lem:2}
The set $\mathcal{Q}$ is non-empty.
\end{lem}
\textit{Proof:} This is immediate from Lemma~\ref{lem:P} and
\ref{lem:R}. We omit the proof for brevity. $\Box$

\begin{prop}
\label{prop:monotonicity}
The function $h^*$ in (\ref{eq:13}) is monotonically increasing in the limit $t\to\infty$.
\end{prop}
\textit{Proof:} From Lemmas~\ref{lem:Q},\ref{lem:P},\ref{lem:R}
and \ref{lem:2}, we have that there exists a $a_2\in\mathcal{Q}$
such that $h(a_2,r)$ is a solution of (\ref{eq:new}) subject to
the boundary conditions (\ref{eq:new2}). From
Proposition~\ref{prop:uniqueness}, we have that (\ref{eq:new}) and
(\ref{eq:new2}) admit a unique solution $h^*$ in the limit
$t\to\infty$. Hence, we deduce that the corresponding shooting
parameter $a^*\in\mathcal{Q}$ i.e. $h^*$ is monotonically
increasing everywhere away from the origin. An immediate consequence of this global monotonicity is $0< h^*(r)<1$ for $r\in (0,\infty)$. $\Box$

\section{The $L\to 0$ limit}
\label{sec:ns}

Consider the Landau-de Gennes energy functional in (\ref{eq:1})
$$\Ical[\Q] = \int_{B(0,R)} \frac{L}{2} |\grad \Q|^2 + f_B(\Q)~dV$$
in the limit $L\to 0^+$. Let $\Q^L$ be an arbitrary solution of the corresponding Euler-Lagrange equations
$$ L \Delta \Q_{ij} = \frac{\partial f_B}{\partial \Q_{ij}} - \frac{1}{3}\frac{\partial f_B}{\partial \Q_{kk}}\delta_{ij} \quad i,j,k=1\ldots 3$$
(where $\frac{1}{3}\frac{\partial f_B}{\partial \Q_{kk}}\delta_{ij}$ is a Lagrange multiplier accounting for tracelessness)
that satisfies an energy bound of the form
\begin{equation}
\label{eq:L1} \Ical[\Q^L] \leq L C(a^2,b^2,c^2,\Omega) \quad \forall L>0
\end{equation}
where $C$ does not depend on $L$ in the limit $L\to 0^+$. Examples of such solutions include global energy minimizers and the radial-hedgehog solution (see (\ref{eq:14})). Then it is intuitively clear that as $L\to 0^+$, such solutions will be \emph{almost} like the bulk energy minimizers in (\ref{eq:3}) and (\ref{eq:4}) (since the energy bound implies $\int_{B(0,R)} f_B(\Q)~dV \leq L C(a^2,b^2,c^2,\Omega)$ as $L\to 0^+$), with the elastic energy density being dominant in the vicnity of defects and interfaces. This has been rigorously established for global energy minimizers in \cite{amaz} and for solutions satisfying the energy bound in (\ref{eq:L1}), in \cite{maj1}.

Consider the ordinary differential equation for $h^*$ in (\ref{eq:11}) and the boundary conditions (\ref{eq:12}) in the limit $L\to 0^+$ i.e.
\begin{eqnarray}
\label{eq:L2}
&& \frac{d^2 h}{dr^2} + \frac{2}{r}\frac{dh}{dr} - \frac{6h}{r^2} = -h + h^3 + \frac{3 h_+}{t}\left(h^3 - h^2\right) \nonumber \\&&
h(0) = 0~ \quad h(r)\to 1~as~r\to\infty
\end{eqnarray}
since the ball radius is inversely proportional to the correlation
length $\xi$ defined in (\ref{eq:nondim}). There is an important
difference between (\ref{eq:L2}) and (\ref{eq:new}) which focusses
on the limit $t\to\infty$. In (\ref{eq:L2}), the term $\frac{3
h_+}{t}\left(h^3 - h^2\right)$ is not necessarily much smaller
than the term $-h + h^3$ on the right-hand side of (\ref{eq:L2}).
Hence, (\ref{eq:L2}) does not have the Ginzburg-Landau structure
as in (\ref{eq:lowtemp2}) and we do not have analogous uniqueness
and global monotonicity results for the radial-hedgehog solution
in the $L\to 0$ limit. We can, however, utilise Ginzburg-Landau
techniques to understand the far-field properties, away from the
isotropic defect core, in the $L\to 0^+$ limit.

We first recall an important result from \cite{amaz} and \cite{maj1} regarding the uniform convergence of $h^*$ everywhere away from the isotropic core.

\begin{prop}\label{prop:mon}\cite{amaz,maj1} For fixed $t>1$ and $L$ sufficiently small,
 there exists $R>0$ such that $h^*(r)\to 1$ uniformly for all $r\geq R$.
\end{prop}

\begin{prop}
\label{prop:farfield} Let $h^*$ be a global minimizer of $I$ in
(\ref{eq:10}), in the class $\Acal_h$, in the limit $L\to 0$. Then
there exists $R_1>0$ such that $h^*$ is monotonically increasing
for all $r\geq R_1$.
\end{prop}

\textit{Proof:} From Proposition~\ref{prop:mon}, we have that
there exists $R_0>0$ such that $h(r)>\frac{1}{2}$ for $r\geq R_0$.
Consider the right-hand side of (\ref{eq:L2}) and define
\begin{equation}
\label{eq:F}
F(h) = h^2 -1 + \frac{3 h_+}{t}(h^2 - h).
\end{equation} Then $F(1)=0$ and $F^{'}(h)>0$ for $h>\frac{1}{2}$.

We prove Proposition~\ref{prop:farfield} by contradiction. We assume that there exists $r_0 >R$, where $R$ is defined in Proposition~\ref{prop:mon}, such that
$$\frac{ d h^*}{ dr}|_{r=r_0} = 0.$$ There are three possibilities for $\frac{d^2 h^*}{d r^2}|_{r=r_0}$ i.e.
(a)$ \frac{d^2 h^*}{d r^2}|_{r=r_0} = 0$, (b) $\frac{d^2 h^*}{d r^2}|_{r=r_0} <0$, and (c) $\frac{d^2 h^*}{d r^2}|_{r=r_0}>0$.

Consider case (a). Then we have from (\ref{eq:L2}) that
\begin{equation}\label{eq:L3}
 \frac{d^2 h^*}{dr^2} + \frac{2}{r}\frac{dh^*}{dr} = h^*\left[ F(h^*) + \frac{6}{r^2}\right] = 0~ at~ r=r_0.
 \end{equation}
 Secondly,
 $$ \frac{d}{dr}\left[ F(h^*) + \frac{6}{r^2}\right] < 0~ at~ r=r_0,$$
 from which we deduce that
 \begin{equation}
 \label{eq:L4}
 F(h^*) + \frac{6}{r^2} > 0 \quad r\in\left(r_0 - \delta, r_0\right)
 \end{equation} for some $\delta>0$.
 We deduce from (\ref{eq:L2}) that
 $$ \frac{d}{dr}\left( r^2 \frac{d h^*}{dr}\right)>0  \quad r\in\left(r_0 - \delta, r_0\right)$$
 so that
 $$ r_0^2 \frac{dh^*}{dr}|_{r=r_0} > (r_0 - \delta)^2 \frac{d h^*}{dr}|_{r_0 - \delta}.$$
 Since $\frac{ d h^*}{ dr}|_{r=r_0} = 0$, we deduce that $\frac{d h^*}{dr}|_{r_0 - \delta} < 0$.
 This necessarily means that there exists a local minimum at $r=r_1 > r_0$,
 since $0\leq h^*(r)\leq 1$ and $h^* \to 1$ as $r\to\infty$.
 We, therefore, have $$\frac{d^2 h^*}{d r^2}|_{r=r_1} > 0$$ or equivalently
 $$ F(h^*(r_1)) + \frac{6}{r_1^2} >0 .$$
 But $$F(h^*(r_1)) + \frac{6}{r_1^2} < F(h^*(r_0)) + \frac{6}{r_0^2} = 0$$
since  $F^{'}(h)>0$ for $h>\frac{1}{2}$ and $h^*(r_1) < h^*(r_0)$. This gives a contradiction and we deduce that $ \frac{d^2 h^*}{d r^2}|_{r=r_0}\neq 0$.

Case (b): We assume that $\frac{d^2 h^*}{d r^2}|_{r=r_0} <0$ and $h^*(r_0)>\frac{1}{2}$ i.e. we have a local maximum at $r=r_0$. The local maximum must be followed by a local minimum at $r=r_1>r_0$, since $0\leq h^*(r)\leq 1 \quad \forall r>0$ and $h^* \to 1$ as $r\to\infty$. Thus,
$$ F(h^*(r_1)) + \frac{6}{r_1^2} >0 $$ by definition of a local minimum from (\ref{eq:L3}).
However $h(r_1)<h (r_0)$ and
$$ F(h^*(r_1)) + \frac{6}{r_1^2} < F(h^*(r_0)) + \frac{6}{r_0^2} < 0$$
yielding a contradiction.

Case (c): We assume that $\frac{d^2 h^*}{d r^2}|_{r=r_0} >0$. Then $\frac{d h^*}{dr}>0$ for $r>r_0>R$ where $R$ has been defined in Proposition~\ref{prop:mon}, since the previous arguments show that we cannot have a point of inflection or a local maximum for $r\geq R$. Then we set $R_1$ in Proposition~\ref{prop:farfield} to be $R_1=r_0$. Proposition~\ref{prop:farfield} now follows. $\Box$

We use the far-field monotonicity established in
Proposition~\ref{prop:farfield} to derive an explicit far-field
expansion for $h^*$ as $r\to\infty$. This expansion is valid in
both the $L\to 0$ and $t\to\infty$ limits.

\begin{prop}
\label{prop:noniso1} Let $h^*$ be a minimizer of $I$ in
(\ref{eq:10}) in the space $\Acal_h$, for a fixed $t>1$, in the
limit $L\to 0$. Then $h^*$ is a non-negative solution of the
following singular ordinary differential equation
\begin{equation}
\label{eq:ns1}
\frac{d^2 h}{dr^2} + \frac{2}{r}\frac{dh}{dr} - \frac{6h}{r^2} = h\left( h^2 - 1 + \frac{3 h_+}{t}\left(h^2 - h\right) \right)
\end{equation}
subject to the boundary conditions
\begin{equation}
\label{eq:ns2}
h(0) = 0 \quad \textrm{and} ~ h(r)\to 1 ~as~r\to\infty.
\end{equation}
We have the following far-field estimates
\begin{equation}
\label{eq:ns2a}
r^2 \left|\frac{d^2 h^*}{dr^2}\right| + r\left|\frac{dh^*}{dr}\right| + \left| 6  - r^2 h^* (1 - h^*)\left( 1 + \left( 1 + \frac{3 h_+}{t}\right) h^* \right)\right| = o(1) \quad r\to\infty.
\end{equation}
\end{prop}


\textit{Proof:}  The proof of Proposition~\ref{prop:noniso1}
follows some of the methods described in a recent paper
\cite{millot} on Ginzburg-Landau theory for three-dimensional
domains.

The bounds (\ref{eq:bounds}) are valid in the $L\to 0$ limit. In fact, they are valid in any limit which translates to an unbounded domain in terms of the dimensionless variables in (\ref{eq:nondim}). In particular, they imply that
$$ 1- \frac{\alpha}{r^2}\leq h^*(r)\leq 1 - \frac{\beta}{r^2}$$
as $r\to\infty$, as shown in (\ref{eq:bc2}), for positive constants $\alpha,\beta$ independent of $L$. As demonstrated in (\ref{eq:far}), this implies
\begin{equation}
\label{eq:r2}
\lim_{r\to\infty} r^2 \frac{ d h^*}{dr} = 0
\end{equation}
and hence
\begin{equation}
\label{eq:r3}
\lim_{r\to\infty} r \frac{ d h^*}{dr} = 0.
\end{equation}

For any $k\in(0,1)$ fixed, we multiply (\ref{eq:ns1}) by $r^2$,
 average over $\left(kR_*, R_* \right)$, take the limit
$R_*\to\infty$ and obtain
\begin{equation}
\label{eq:ns6a} \frac{1}{(1-k)R_*}\int_{k
R_*}^{R_*}\frac{d}{dr}\left(r^2 \frac{dh^*}{dr}\right) dr +
\frac{1}{(1-k) R_*}\int_{k R_*}^{R_*} r^2 h^*(r)\left( 1 -
h^*(r)\right)\left( 1 + \left(1 + \frac{3 h_+}{t}\right)h^*\right)
dr= \frac{6}{(1-k)R_*}\int_{k R_*}^{R_*} h^*(r)~dr.
\end{equation}

In the limit $L\to 0^+$, $h^* \to 1$ uniformly as $r\to\infty$ from Proposition~\ref{prop:mon} and using (\ref{eq:r2}), we obtain the following sequence of inequalities
\begin{equation}
\label{eq:ns7}
\limsup_{R_*\to\infty}k^2 R_*^2 \left(1 - h^*(R_*))\right)\left( 1 + \left(1 + \frac{3 h_+}{t}\right)h^*(R_*)\right) \leq 6 \leq \liminf_{R_*\to\infty}R_{*}^2 \left(1 - h^*(k R_*)\right)\left( 1 + \left(1 + \frac{3 h_+}{t}\right)h^*(k R_*)\right).
\end{equation}
It immediately follows that
\begin{equation}
\label{eq:ns8} r^2\left(1 - h^*(r))\right)\left( 1 + \left(1 +
\frac{3 h_+}{t}\right)h^*(r)\right) \to 6
\end{equation}
uniformly in the limit $r\to\infty$.

Finally, using the estimates (\ref{eq:r3}) and (\ref{eq:ns8}) in (\ref{eq:ns1}), we deduce that
 \begin{equation}
\label{eq:ns9}
r^2 \left|\frac{ d^2 h^*}{dr^2}\right| \to 0
\end{equation}
uniformly in the limit $r \to\infty$.
Proposition~\ref{prop:noniso1} now follows. $\Box$

One immediate consequence of (\ref{eq:ns2a}) is that
\begin{equation}
\label{eq:ns10}
h^*(r) = 1 - \frac{6}{r^2\left(2 + \frac{3 h_+}{t}\right)} + o\left(\frac{1}{r^2}\right) \quad r\to\infty.
\end{equation}
Although this information is qualitatively contained in
(\ref{eq:bc2}), (\ref{eq:ns10}) is a stronger result since it is
an exact expression that captures the effects of geometry and the
temperature on the far-field structure. Further, (\ref{eq:ns10})
yields estimates for the higher-order derivatives of $h^*$ as
$r\to\infty$ and this information cannot be immediately inferred
from (\ref{eq:bc2}).

\begin{prop}
\label{prop:stability1} Let $\Q^L$ denote the radial-hedgehog
solution in (\ref{eq:13}) for a fixed $L>0$ and fixed $t>1$. Then
$\left\{\Q^L\right\} \to \Q^0$ in
$W^{1,2}\left(B(0,R);\bar{S}\right)$ as $L\to 0^+$, where $\Q^0
=\sqrt{\frac{3}{2}} \left(\frac{\rvec\otimes \rvec}{r^2} -
\frac{1}{3}\I \right)$ \cite{amaz}. As $L\to 0^+$, the
radial-hedgehog solution $\Q^L$ is stable against all small
far-field perturbations $\Pvec$ that satisfy
\begin{eqnarray}
\label{eq:s1} && \Pvec(\rvec) = 0 \quad \rvec\in
B(0,R_L),\nonumber\\
&& \Pvec(\rvec) = 0 \quad \rvec\in\partial B(0,R)
\end{eqnarray}
and $R_L$ is sufficiently large. In other words, $\Q^L$ is locally
stable against perturbations which are localized outside the
isotropic core around the origin, in the limit $L\to 0^+$.
\end{prop}

\textit{Proof:} Consider the dimensionless free energy in
(\ref{eq:nondim2})
$$ I[\Q] = \int_{B(0,R)} \frac{1}{2}|\grad \Q|^2 + f_B(\Q)~dV$$
where $$ f_B = -\frac{1}{2}\textrm{tr}\Q^2 -
\frac{\sqrt{6}h_+}{t}\textrm{tr}\Q^3 +
\frac{h_+^2}{2t}\left(\textrm{tr}\Q^2\right)^2 + C(t) \geq 0$$ and
$$ f_B(\Q) = 0 \Longleftrightarrow \Q =
\sqrt{\frac{3}{2}}\left(\nvec\otimes\nvec -
\frac{\mathbf{I}}{3}\right).$$

We consider a small perturbation $\Pvec$ that satisfies
(\ref{eq:s1}). Define
\begin{equation}
\label{eq:p1} \Q^{\eps}_{ij} = \Q^L_{ij} + \eps\Pvec_{ij} \quad
0<\eps<<1.
\end{equation}
We compute the second variation $$\frac{d^2 I[\Q^{\eps}]}{d
\eps^2}|_{\eps=0}.$$ A direct computation shows that
\begin{equation}
\label{eq:p2} \frac{d^2 I[\Q^{\eps}]}{d \eps^2}|_{\eps=0} =
\int_{B(0,R)\setminus B(0,R_L)} \left|\grad \Pvec\right|^2 +
\frac{\partial^2 f_B}{\partial \Q^L_{ij} \partial
\Q^{L}_{pq}}\Pvec_{ij}\Pvec_{pq}~dV
\end{equation}
since $\Pvec=0$ on $B(0,R_L)$.

From Proposition~\ref{prop:noniso1}, we have that as $L\to 0$,
\begin{equation}
\label{eq:p3} \Q^L (\rvec) = \sqrt{\frac{3}{2}}\left(1  -
\frac{\gamma}{r^2} +
o(\frac{1}{r^2})\right)\left(\frac{\rvec}{r}\otimes
\frac{\rvec}{r} - \frac{\mathbf{I}}{3}\right) \quad r\to\infty.
\end{equation} Equation (\ref{eq:p3}) can be written as
\begin{equation}
\label{eq:p4} \Q^L (\rvec) = \Q^0(\rvec)
-\frac{\gamma}{r^2}\left(\frac{\rvec}{r}\otimes \frac{\rvec}{r} -
\frac{\mathbf{I}}{3}\right) + o(\frac{1}{r^2}) \quad r\to\infty
\end{equation}
where $\Q^0$ is a bulk energy minimizer by definition (see
preceding comments).

We perform a Taylor expansion of $\frac{\partial^2 f_B}{\partial
\Q^L_{ij} \partial \Q^{L}_{pq}}$ around $\Q^0$ to obtain
\begin{equation}
\label{eq:p5} \frac{\partial^2 f_B}{\partial \Q^L_{ij} \partial
\Q^{L}_{pq}}\Pvec_{ij}\Pvec_{pq} =\frac{\partial^2 f_B}{\partial
\Q^0_{ij} \partial \Q^{0}_{pq}}\Pvec_{ij}\Pvec_{pq} -
\frac{\gamma}{r^2}
\left(\frac{\rvec_{\alpha}\rvec_{\beta}}{r^2}-\frac{\delta_{\alpha\beta}}{3}\right)
\frac{\partial^3 f_B}{\partial \Q^0_{\alpha\beta}\partial
\Q^0_{ij}
\partial \Q^{0}_{pq}} + o\left(\frac{1}{r^2}\right) \quad
r\to\infty.
\end{equation}
Finally, we note that
$$\frac{\partial^2 f_B}{\partial
\Q^0_{ij} \partial \Q^{0}_{pq}}\Pvec_{ij}\Pvec_{pq} > \eta(t)$$
where $\eta$ is independent of $r$ (since $\Q^0$ is a global
minimizer of $f_B$ by definition) and
$$\left|\frac{\gamma}{r^2}
\left(\frac{\rvec_{\alpha}\rvec_{\beta}}{r^2}-\frac{\delta_{\alpha\beta}}{3}\right)
\frac{\partial^3 f_B}{\partial \Q^0_{\alpha\beta}\partial
\Q^0_{ij}
\partial \Q^{0}_{pq}} \right|\leq \frac{\psi(t)}{R_L^2} \quad
r\geq R_L$$ since the derivatives of $f_B$ can be bounded
independently of $R_L$. Combining the above, we have
\begin{equation}
\label{eq:p6} \frac{\partial^2 f_B}{\partial \Q^L_{ij} \partial
\Q^{L}_{pq}}\Pvec_{ij}\Pvec_{pq} \geq \eta(t) -
\frac{\psi^{'}(t)}{R_L^2} >0
\end{equation}
for a fixed $t>1$ and $R_L$ sufficiently large. Substituting
(\ref{eq:p6}) into (\ref{eq:p2}), we deduce that
\begin{equation}
\label{eq:p7} \frac{d^2 I[\Q^{\eps}]}{d \eps^2}|_{\eps=0} > 0
\end{equation}
for perturbations $\Pvec$ satisfying (\ref{eq:s1}). The positivity
of the second variation ensures that the radial-hedgehog solution
is locally stable against perturbations $\Pvec$ satisfying
(\ref{eq:s1}), in the limit $L \to 0^+$.
Proposition~\ref{prop:stability1} now follows. $\Box$

\section{A general stability result}
\label{sec:stability}

In the previous sections, we have demonstrated that the
radial-hedgehog solution is unstable with respect to biaxial
perturbations, localised around the isotropic core, for $R$ and
$t$ sufficiently large. As already stated, the $R\to\infty$ limit
is equivalent to either the $L\to 0$ or $t\to\infty$ limits or
both. We have also shown that the radial-hedgehog solution is
locally stable with respect to far-field perturbations in the $R
\to\infty$ limit. We conclude by deriving a general local
stability result that is not restricted to the limits $t\to\infty$
or $L \to 0$. We note that Proposition~\ref{prop:stability2} is known from numerical investigations (see \cite{sonnet,mkaddemgratland1, mkaddemgartland2}) and we present a proof partly for completeness and partly this proof gives greater insight into how the elastic constant, temperature and ball radius collectively 
affect stability properties.
\begin{prop}
\label{prop:stability2} Let $B(0,R)$ denote a ball of radius $R$
centered at the origin in $\mathbb{R}^3$. The corresponding
radial-hedgehog solution $\Q^*_R$ is stable against all small,
smooth perturbations of the form
\begin{equation}
\Q = \Q^*_R + \eps \Pvec \label{eq:s13}
\end{equation}
where $\eps\in \Rr$, $|\eps|<<1$, $\Pvec\in \bar{S}$ and $\Pvec=0$
on $\partial B(0,R)$, provided that the radius $R$ is sufficiently
small i.e.
\begin{equation}
\label{eq:s12} R^2 < \frac{1}{4}\left(\frac{1}{1 +
\frac{4\sqrt{6}h_+}{t}}\right).
\end{equation}
In terms of the original variables defined in (\ref{eq:nondim}),
(\ref{eq:s12}) is equivalent to
\begin{equation}
\label{eq:s13b} R_{real}^2 < \frac{\xi^2}{ 4 t}\left( \frac{1}{1 +
\frac{4\sqrt{6}h_+}{t}}\right)
\end{equation}
where $\xi$ is the correlation length defined in
(\ref{eq:nondim}).
\end{prop}

\textit{Proof:} The results in Proposition~\ref{prop:1} are true
for any $R>0$ i.e. for every $R>0$, we are guaranteed the
existence of a radial-hedgehog solution $\Q^*_R$ of the form
(\ref{eq:13}), that satisfies the energy bound (\ref{eq:14}) and
the inequalities (\ref{eq:15}). Consider the dimensionless free
energy in (\ref{eq:nondim4}) and introduce the change of variable
$$\hat{r} = \frac{r}{R}$$
so that the free energy becomes
\begin{equation}
\label{eq:s14} I[\Q] =
\int_{0}^{2\pi}\int_{0}^{\pi}\int_{0}^{1}\left\{\frac{1}{2}|\grad
\Q|^2 - \frac{R^2}{2}\textrm{tr}\Q^2 -
\frac{\sqrt{6}h_+}{t}R^2\textrm{tr}\Q^3 +
\frac{h_+^2}{2t}R^2\left(\textrm{tr}\Q\right)^2 + R^2 C(t)\right\}
\hat{r}^2 \sin\theta d\hat{r}d\theta d\phi.
\end{equation}

We consider small perturbations
\begin{equation} \label{eq:s15}
\Q_{\eps} = \Q^*_R +\eps \Pvec \quad 0<\eps<<1
\end{equation} such that $\Pvec=0$ on $\partial B(0,R)$. We
compute the second variation of the Landau-de Gennes energy
functional
\begin{equation}
\label{eq:s16} \frac{d^2}{d\eps^2}I[\Q_\eps]|_{\eps=0} =
\int_{0}^{2\pi}\int_{0}^{\pi}\int_{0}^{1}\left\{ |\grad \Pvec|^2 -
R^2 |\Pvec|^2 - \frac{6\sqrt{6}h_+}{t}R^2
\Pvec_{ij}\Pvec_{jp}\Q^*_{R_{pi}} + \frac{h_+^2
R^2}{2t}\left[8\left(\Q^*_R\cdot\Pvec\right)^2 +
4|\Pvec|^2|\Q^*_R|^2\right]\right\}~dV
\end{equation} where $dV =\hat{r}^2 \sin\theta d\hat{r}d\theta
d\phi$.

We, next, make an elementary observation
$$ \Pvec_{ij}\Pvec_{jp}\Q^*_{R_{pi}} = h^*(r)\left[\rvec_i
\Pvec_{ij}\rvec_p\Pvec_{pj}/r^2 - \left|\Pvec\right|^2/3\right]
\leq \frac{2}{3}|\Pvec|^2 $$ so that
\begin{equation}
\label{eq:s17} \frac{d^2}{d\eps^2}I[\Q_\eps]|_{\eps=0} \geq
\int_{0}^{2\pi}\int_{0}^{\pi}\int_{0}^{1}\left\{ |\grad \Pvec|^2
\hat{r}^2 - R^2 \hat{r}^2|\Pvec|^2 -
\frac{4\sqrt{6}h_+}{t}R^2\hat{r}^2 |\Pvec|^2
\right\}\sin\theta~d\hat{r}~d\theta~d\phi.
\end{equation}
We note that
$$|\Pvec|^2 \geq \left(\frac{\partial \Pvec}{\partial \hat{r}}\right)^2$$ and
use the following inequality from \cite{chipark,kinderlehrer}
$$ \int_{0}^{1}\tau^2\left(\frac{\partial\alpha}{\partial
\tau}\right)^2~d\tau \geq
\frac{1}{4}\int_{0}^{1}\alpha^2(\tau)~d\tau$$ for a real-valued
function $\alpha$ defined on the interval $[0,1]$. Substituting
the above inequality in (\ref{eq:s17}), we have that
\begin{equation}
\label{eq:s17} \frac{d^2}{d\eps^2}I[\Q_\eps]|_{\eps=0} \geq
\int_{0}^{2\pi}\int_{0}^{\pi}\int_{0}^{1}\left\{
\frac{1}{4}|\Pvec|^2 - |\Pvec|^2 R^2\left(1 +
\frac{4\sqrt{6}h_+}{t}\right)
\right\}\sin\theta~d\hat{r}~d\theta~d\phi
\end{equation} since $\hat{r}\leq 1$. It follows that
$$\frac{d^2}{d\eps^2}I[\Q_\eps]|_{\eps=0} > 0$$
if
\begin{equation}
\label{eq:s18} R^2 < \frac{1}{4}\frac{1}{1 +
\frac{4\sqrt{6}h_+}{t}}
\end{equation}
or equivalently if
\begin{equation}
\label{eq:s19} R_{real}^2 < \frac{\xi^2}{4t}\left(\frac{1}{1 +
\frac{4\sqrt{6}h_+}{t}}\right)
\end{equation}
where $R_{real} = \frac{\xi}{\sqrt{t}}R$ from (\ref{eq:nondim})
and (\ref{eq:nondimnew}). Proposition~\ref{prop:stability2} now
follows. $\Box$

\section{Discussion}
\label{sec:discussion}

This paper aims to build a self-contained and rigorous mathematical framework for the study of the radial-hedgehog solution within the Landau-de Gennes theory for nematic liquid crystals and elucidate the analogies between the mathematical formulation of defects in the Landau-de Gennes framework and defects in the Ginzburg-Landau theory of superconductivity. These analogies need to be highlighted in the applied mathematics literature, so that mathematical techniques from other branches of condensed matter science can be effectively used in the context of liquid crystals. We study radial-hedgehog solutions on spherical droplets subject to \emph{homeotropic} anchoring or \emph{strong radial} anchoring conditions and define a radial-hedgehog solution to be an energy minimizer within the class of spherically symmetric uniaxial solutions as demonstrated in Proposition~\ref{prop:1}. We consider two different limits in this paper: (a) the low--temperature limit $t\to\infty$ and (b) the vanishing core limit $L\to 0^+$. For completeness, we summarize the validity of the different results in this paper in different parameter regimes. Proposition~\ref{prop:1} and Proposition~\ref{prop:iso2} are valid in all parameter regimes i.e. we are always guaranteed the existence of a radial-hedgehog solution that satisfies the energy bound (\ref{eq:14}) and whose scalar order parameter is constrained by the inequalities (\ref{eq:15}). We have a single isolated isotropic point at the origin by definition and we always have a series expansion near the origin that involves only even powers of $r$ as $r\to 0$. However, the bounds (\ref{eq:bounds}) are only valid in the limit $R\to\infty$, where $R$ is the re-scaled ball radius. Recalling the definition of the dimensionless variables in (\ref{eq:nondim}) and (\ref{eq:nondimnew}), the limit $R\to\infty$ is equivalent to either the low--temperature limit or the vanishing core limit. Propositions~\ref{prop:noniso1} and Proposition~\ref{prop:stability1} are valid in the $L\to 0$ limit, whereby it is difficult to establish rigorous results about the defect core but the governing equation has a Ginzburg-Landau structure away from the origin. This Ginzburg-Landau structure gives us a good grip on the far-field properties i.e. uniform convergence of the scalar order parameter away from the origin, far-field monotonicity and explicit far-field expansions for the scalar order parameter. The limit $t\to\infty$ has a Ginzburg-Landau structure (see (\ref{eq:lowtemp}) and (\ref{eq:lowtemp2})) and we can exploit Ginzburg-Landau techniques to prove global properties. Propositions~\ref{prop:uniqueness}, \ref{prop:monotonicity}, \ref{prop:noniso1}, \ref{prop:stability1} and \ref{prop:biaxial} hold in the $t\to\infty$ limit and we demonstrate the manifestation of biaxial instabilities localised near the isotropic core in this regime. Proposition~\ref{prop:uniqueness} is an example of how Ginzburg-Landau techniques can be used to prove results on multiplicity of solutions and in Proposition~\ref{prop:monotonicity}, we appeal to shooting arguments  which have not been used previously in the Landau-de Gennes context. Proposition~\ref{prop:stability2} is a general result that identifies a relationship between the elastic constant $L$, the reduced temperature $t$ and the ball radius $R$ that guarantees local stability of the radial-hedgehog solution against all perturbations.

In \cite{brezis2}, H. Brezis postulated the following problem in the context of Ginzburg-Landau theory for superconductors:
for maps $\uvec:\Rr^3\to \Rr^3$, is any solution of the system
\begin{equation}
\label{eq:dis1}
\Delta \uvec + \uvec\left(1 - |\uvec|^2\right) = 0\end{equation} satisfying $|\uvec(\rvec)| \to 1 ~\textrm{as}~ |\rvec| → +\infty$ (possibly with a “good” rate of convergence) and ${deg}_\infty \uvec = ±1$ of the form
\begin{equation}
\label{eq:dis2}\mathbf{U}(\rvec) = \frac{\rvec}{r}f(r) \end{equation}
for a unique function $f$ vanishing at zero and increasing to one at infinity. In \cite{millot}, the authors show that every non-constant local minimizer of the Ginzburg-Landau energy functional associated with (\ref{eq:dis1}),
$$ E(\uvec,\Omega) :=\int_{\Omega}\frac{1}{2}|\grad \uvec|^2 + \frac{1}{4}\left(1 - |\uvec|^2\right)^2~dV$$ is of the form (\ref{eq:dis2}), up to a translation on the domain and an orthogonal transformation on the image.
For nematic liquid crystals, the corresponding problem translates to: is any uniaxial solution of (\ref{eq:7}) necessarily of the form (\ref{eq:13}) i.e. are radial-hedgehog solutions the only possible uniaxial solutions of the system (\ref{eq:7}) in $\Rr^3$? If so, then we will have a complete characterization of all admissible uniaxial solutions and the interplay between uniaxiality and biaxiality can be partially understood in terms of the comparitively tractable radial-hedgehog problem. We expect that the methods in \cite{millot} will not readily transfer to the Landau-de Gennes framework and there will be analogies only in certain parameter regimes, such as the $t\to\infty$ limit studied in this paper.

Finally, we compare our results with previous work in this area. In \cite{sonnet}, the authors carry out detailed numerical investigations of equilibrium configurations within spherical droplets subject to strong radial anchoring conditions and find that the radial-hedgehog solution only occurs either in very small droplets or very close to the nematic-isotropic transition temperature; the symmetry-breaking biaxial torus solution is energetically preferable everywhere else. This is consistent with Proposition~\ref{prop:stability2} where we demonstrate local stability of the radial-hedgehog solutions for droplets with radius comparable to the nematic correlation length $\xi$. This is also consistent with Proposition~\ref{prop:biaxial} where we demonstrate that the radial-hedgehog solution cannot be a global energy minimizer for large droplets in the low-temperature limit. In \cite{rossovirga}, the authors work within the Lyuksyutov constraint, which requires that $\textrm{tr}\Q^2(\rvec) = \frac{2}{3}s_+^2$ for $\rvec\in B(0,R)$, where $s_+$ has been defined in (\ref{eq:4}). They demonstrate that the radial-hedhegog solution is always locally \emph{unstable} within the one-constant approximation for the elastic energy density i.e. when the elastic energy density is simply taken to $|\grad \Q|^2$, as has been done in this paper. This resuly is evidently in agreement with Proposition~\ref{prop:biaxial} and does \emph{not} contradict Proposition~\ref{prop:stability2} where we demonstrate local stability in balls of sufficiently small radius. The Lyuksyutov constraint is valid in the $R\to\infty$ limit or for balls of sufficiently large radius and hence, Proposition~\ref{prop:stability2} is outside the remit of this instability result. The analogies of this work with the results reported in \cite{mkaddem&gartland1, mkaddem&gartland2} have already been mentioned. The authors numerically study the stability of the radial-hedgehog solution as a function of the ball radius, reduced temperature and elastic constants in \cite{mkaddem&gartland2,mkaddem&gartland1} ( along with analysis of the biaxial instabilites as has been mentioned in Section~\ref{sec:isotropic}) and in the one-constant case, our results are qualitatively in agreement with the phase diagrams in \cite{mkaddem&gartland2, mkaddem&gartland1}. It would be interesting to see if (\ref{eq:s19}) can yield a qualitative fit to the region of local stability obtained in \cite{mkaddem&gartland2, mkaddem&gartland1}. While careful attention is paid to the effect of elastic constants in some of the previous work, we focus on the one-constant case. This is primarily because the one-constant case has a much more tractable mathematical structure than the unequal constant case and is the best paradigm for illustrating the efficacy of Ginzburg-Landau techniques in the Landau-de Gennes framework. The unequal elastic constant case will be considered in future work.

\section*{Acknowledgments} This publication is based on work supported by
Award No. KUK-C1-013-04 , made by King Abdullah University of
Science and Technology (KAUST) to the Oxford Centre for
Collaborative Applied Mathematics. The author gratefully
acknowledges stimulating discussions with Maria Aguareles, Chong Luo, Luc Nguyen and Arghir Zarnescu. We thank Luc Ngyuen and Arghir Zarnescu for helpful comments and suggestions regarding Proposition~\ref{prop:1}.
\begin{small}

\end{small}
\end{document}